\documentclass[11pt,a4paper]{article}
\usepackage[T1]{fontenc}
\usepackage{lmodern}
\usepackage[margin=1in]{geometry}
\usepackage{amsmath,amssymb,amsthm}
\usepackage{enumitem}
\usepackage{tikz}
\usetikzlibrary{arrows.meta}
\usepackage[colorlinks=true,linkcolor=blue,citecolor=blue,urlcolor=blue]{hyperref}

\theoremstyle{plain}
\newtheorem{theorem}{Theorem}[section]
\newtheorem{lemma}[theorem]{Lemma}
\newtheorem{proposition}[theorem]{Proposition}
\newtheorem{corollary}[theorem]{Corollary}
\newtheorem{conjecture}[theorem]{Conjecture}

\newcommand{\odiam}{\overrightarrow{\mathrm{diam}}}
\DeclareMathOperator{\dist}{dist}
\DeclareMathOperator{\diam}{diam}
\newcommand{\Nbad}{N_{\mathrm{bad}}}

\title{An Erd\H{o}s–Pach–Pollack–Tuza theorem for oriented diameter}
\author{Jiangdong Ai\thanks{School of Mathematical Sciences and LPMC, Nankai University. \texttt{jd@nankai.edu.cn}. Funded by the National Natural Science Foundation of China (No.~12522117) and the Fundamental and Interdisciplinary Disciplines Breakthrough Plan of the Ministry of Education of China (JYB2025XDXM207).}
\and Bin Chen\thanks{School of Mathematics and Statistics, Fuzhou University, Fujian, China. \texttt{cbfzu03@163.com}. 
Funded by the National Natural Science Foundation of China (Nos.~12471336,~12501473).}
\and Yaokun Feng\thanks{Center for Combinatorics and LPMC, Nankai University, and Department of Applied Mathematics, The Hong Kong Polytechnic University, Hong Kong, China. \texttt{1031679547@qq.com}.}}
\date{}

\begin{document}

\maketitle

\begin{abstract}
The oriented diameter of a bridgeless graph is the minimum diameter of its strong orientations. Erd\H{o}s, Pach, Pollack and Tuza proved that every connected graph of order $n$ and minimum degree $\delta\ge 2$ has diameter at most $3n/(\delta+1)-1$, which is sharp up to an additive constant. Cochran and Wang conjectured that the oriented diameter obeys the same bound up to an absolute additive constant. We prove this conjecture: every bridgeless graph of order $n$ and minimum degree $\delta$ has a strong orientation of diameter at most $3n/(\delta+1)+10^{17}$. Thus the bound of Erd\H{o}s, Pach, Pollack and Tuza for the diameter carries over to the oriented diameter, with only the additive constant changed; in particular, the coefficient $3$ is best possible.


\end{abstract}

\section{Introduction}

Throughout, all graphs are finite and simple unless stated otherwise; multigraphs arise only from contractions, and they are always loopless. An \emph{orientation} of a graph $G$ replaces each edge by one of its two possible arcs, and it is \emph{strong} if every ordered pair of vertices is joined by a directed path. A \emph{bridge} is an edge
of a connected graph after whose removal the graph is no longer connected, and a graph is said to be \emph{bridgeless} if it is
connected and contains no bridge. Robbins' theorem~\cite{Rob39} states that a connected graph has a strong orientation if and only if it is bridgeless. For a strong orientation $D$, we write $\dist_D(x,y)$ for the length of a shortest directed path from $x$ to $y$, and $\diam(D)$ for the maximum of $\dist_D(x,y)$ over all ordered pairs $(x,y)$. The \emph{oriented diameter} $\odiam(G)$ of a bridgeless graph $G$ is the minimum of $\diam(D)$ over all strong orientations $D$ of $G$.

The study of distances in orientations goes back to Chvátal and Thomassen~\cite{CT78}. Among other results, they proved that every bridgeless graph of radius $r$ has an orientation of radius at most $r^2+r$; this radius theorem is used in the last step of our proof. In this paper we bound the oriented diameter in terms of the order and the minimum degree. For undirected graphs, Erd\H{o}s, Pach, Pollack and Tuza~\cite{EPPT89} showed that every connected graph of order $n$ and minimum degree $\delta\ge 2$ has diameter at most $\lfloor 3n/(\delta+1)\rfloor-1$, and this is sharp up to an additive constant. For orientations, Bau and Dankelmann~\cite{BD15} proved that every bridgeless graph $G$ of order $n$ and minimum degree $\delta$ satisfies
\[
\odiam(G)\le \frac{11n}{\delta+1}+9 .
\]
Surmacs~\cite{Sur17} improved the coefficient $11$ to $7$, and Cochran~\cite{Coc24} obtained a bound of the form $(6+\varepsilon)n/(\delta+1)+C_\varepsilon$ for every fixed $\varepsilon>0$. Recently, Cochran and Wang~\cite{CW26} proved that for every $\varepsilon>0$ there is a constant $C_\varepsilon$ such that every bridgeless graph $G$ of order $n$ and minimum degree $\delta\ge 3$ satisfies $\odiam(G)\le (3+\varepsilon)n/(\delta-2)+C_\varepsilon$, provided that $n$ is sufficiently large. They proposed the following conjecture.

\begin{conjecture}\label{conj:CW}
There is an absolute constant $C$ such that every bridgeless graph $G$ of order $n$ and minimum degree $\delta\ge 2$ satisfies
\[
\odiam(G)\le \frac{3n}{\delta+1}+C .
\]
\end{conjecture}

The coefficient $3$ in Conjecture~\ref{conj:CW} cannot be lowered, as already shown by Bau and Dankelmann~\cite{BD15}. Indeed,  Bau and Dankelmann constructed a family of graphs with minimum degree $\delta\geq 4$ in which every graph $G$ has oriented diameter at least $\frac{3n}{\delta+1}-3$, and  Cochran and Wang~\cite{CW26} gave another family of graphs with minimum degree $\delta\geq 4$ such that every member $G$ satisfies that its oriented diameter is also at least $\frac{3n}{\delta+1}-3$.  Here we give a better construction than the previous ones. Let $G$ be a graph obtained from a cycle $C= v_1 v_2 \ldots v_{3k} v_1$ by replacing $v_{3i}$ with a complete graph $K^{3i}$ on $\delta-1$ vertices, and each vertex of $K^{3i}$ is connected to $v_{3i-1}$ and $v_{3i+1}$, where $k\geq 2$, $1\leq i\leq k$, $\delta\geq 2$, and indices are taken modulo $3k$. Observe that $G$ is a graph of order $n=k(\delta+1)$ vertices and minimum degree $\delta(G)=\delta$ satisfying that
$\overrightarrow{\operatorname{diam}}(G)\geq \frac{3n}{\delta+1}-1$. Indeed, $G$ is bridgeless, and any two of the $k$ edges $v_{3i+1}v_{3i+2}$ form an edge cut, so in every strong orientation of $G$ these edges all point in the same direction around the cycle. If they point from $v_{3i+1}$ to $v_{3i+2}$, then every directed path from $v_{3i+2}$ to $v_{3i+1}$ passes in order through all $k$ cliques $K^{3j}$, using at least two edges at each of them and one edge between consecutive ones, and therefore has length at least $2k+(k-1)=3k-1=\frac{3n}{\delta+1}-1$; the other direction is symmetric.

Our main result confirms Conjecture~\ref{conj:CW}.

\begin{theorem}\label{thm:main}
Every bridgeless graph $G$ of order $n$ and minimum degree $\delta$ satisfies
\[
\odiam(G)\le \frac{3n}{\delta+1}+10^{17}.
\]
\end{theorem}

Thus the bound of Erd\H{o}s, Pach, Pollack and Tuza~\cite{EPPT89} for the diameter carries over to the oriented diameter of bridgeless graphs, with only the additive constant changed. Note that this concerns the extremal bound in terms of $n$ and $\delta$; for an individual graph the oriented diameter can be much larger than the diameter. We have not tried to optimise the additive constant. We keep all constants explicit in order to make clear that the additive term depends neither on $n$ nor on $\delta$; the proof actually gives $3n/(\delta+1)+3.51\cdot 10^{16}$.

We now describe the proof. Write $d=\delta+1$. If $P$ is a shortest path in a graph of minimum degree $\delta$, then every vertex $z$ satisfies $|N[z]\cap V(P)|\le 3$, and counting incidences gives $d|V(P)|\le 3|N[V(P)]|$. This is the source of the coefficient $3$ in~\cite{EPPT89}. In an orientation we need two edge-disjoint paths instead of one. We grow a \emph{core}, that is, a vertex set $K$ together with a strong orientation of the subgraph induced by $K$, and we keep the diameter of this orientation at most $3|N[K]|/d$. To enlarge the core towards a vertex $x$, we contract $K$ to a single vertex and choose two edge-disjoint paths from the contracted vertex to $x$ of minimum total length. Their union is a chain of cycles, and we orient every cycle with the help of short paths between its two branches. The diameter of the enlarged core then increases by at most $3m/d$ plus an absolute constant, where $m$ is the number of new vertices in the closed neighbourhood of the core.

The additive constant is absorbed by a strict gain in the counting. The incidence inequality above is close to equality only when the paths run through long sequences of cliques in which every three consecutive cliques have exactly $d$ vertices. We call such sequences exact strips, and we remove long exact strips in advance by a reduction which changes the oriented diameter by exactly the right amount. Without long exact strips, the neighbourhood count beats the incidence inequality by a fixed proportion of the length of the paths. Consequently, an extension towards a vertex whose two paths have total length greater than $3.6\cdot10^7d$ costs at most $3m/d$.

After this stage, every vertex has two edge-disjoint paths to the contracted core of total length $O(d)$, but the number of remaining vertices may still be large. We then select such pairs of paths in batches. A pair that is still needed after the earlier pairs of the batch have been contracted has many new neighbours, so a batch at length $L$ contains $O(q/(dL))$ pairs, where $q$ is the number of vertices outside the closed neighbourhood of the core. Each pair forms a closed trail through the contracted root, and we orient all trails of a batch at once by random priorities; the diameter of their union is $O(L\log(1+k))$ for $k$ trails. Halving $L$ from stage to stage gives a total cost of $O(d[1+\log(1+q/d^2)])$, which is at most $3q/d$ when $q/d^2$ is large. When $q/d^2$ is bounded, a separate density bound of the form $3n/d+O(1+n/d^2)$, which preserves the core, finishes the orientation with an absolute error. For $4\le\delta\le6$ the same stopping condition already bounds the distance to the core by an absolute constant, and the radius theorem of Chvátal and Thomassen completes the orientation. Minimum degree three requires a separate fragment count before this radius completion.

Two points need care. All paths of a batch are chosen with respect to the same root before any direction is assigned, so we never need individually chosen directions to agree on shared edges. Moreover, each vertex is charged only when it first enters the closed neighbourhood of the growing core. The choice between the two ways of finishing depends only on the graph and not on the ordered pair whose distance is estimated, so a single orientation satisfies all distance bounds.

The paper is organised as follows. Section~\ref{sec:paths} collects properties of pairs of edge-disjoint paths of minimum total length. Section~\ref{sec:strips} removes long exact strips. Section~\ref{sec:counting} proves the strict counting lemma and the extension of an oriented core for $\delta\ge 4$, and Section~\ref{sec:small} treats $\delta=3$ and constructs the initial core. Section~\ref{sec:density} proves the density bound, Section~\ref{sec:trails} orients batches of closed trails through one root, and Section~\ref{sec:completion} completes the proof of Theorem~\ref{thm:main}.

\section{Pairs of edge-disjoint paths}\label{sec:paths}

We first fix notation. For a graph $G$ and a set $X\subseteq V(G)$, we write $N_G[X]$ for the set consisting of $X$ and all vertices adjacent to a vertex of $X$, and $N_G^j[X]$ for the set of vertices at distance at most $j$ from $X$; for a single vertex we omit the braces. Degrees and neighbourhoods always refer to the fixed graph indicated by the subscript, which we omit when it is clear from the context. A path has no repeated vertex, and a trail has no repeated edge but may repeat vertices. The length of a path, trail or walk is its number of edges or arcs. When a connected vertex set $K$ is contracted to a single vertex, parallel edges are kept, each with its original end in $K$, and loops are deleted. Natural logarithms are denoted by $\ln$.

A \emph{core} is a nonempty vertex set $K$ together with a prescribed strong orientation of the subgraph induced by $K$; a single vertex is a core of diameter $0$. When a core is enlarged, all directions prescribed inside the old core are kept. Throughout the paper the minimum-degree parameter $\delta$ is fixed before any reduction, and $d=\delta+1$. In reduced graphs we only require minimum degree at least $\delta=d-1$, and we never replace $d$ by the actual minimum degree.

\begin{lemma}\label{lem:lift}
Let $K$ be a core of a connected graph $H$, and let its prescribed orientation have diameter $\Delta_K$. Every strong orientation $D$ of the loopless multigraph $H/K$, with physical edges recorded at their actual ends in $K$, lifts to a strong orientation of $H$ preserving the prescribed directions and having diameter at most $\Delta_K+\diam(D)$.
\end{lemma}

\begin{proof}
Keep the prescribed orientation inside $K$ and give every other edge the direction of its physical image in $D$. For an ordered pair not wholly in $K$, choose a shortest directed path between its images in $D$. This path visits the contracted root at most once. Replace that visit, if present, by a directed path inside $K$ between the recorded ends, of length at most $\Delta_K$. If an endpoint lies in $K$, use the corresponding initial or final path inside $K$. Pairs wholly in $K$ use their prescribed routes. These choices use one fixed orientation, and give the stated bound.
\end{proof}

We use two edge-disjoint paths of minimum total length rather than two separately chosen paths. Such pairs are standard objects in disjoint-path algorithms; see Suurballe and Tarjan~\cite{ST84}. The optimality conditions for minimum-cost flows~\cite[Chapter~9]{AMO93} give the following description. All edges of a multigraph are regarded as distinct physical edges, including edges with the same ends.

\begin{lemma}\label{lem:potential}
Let $\Gamma$ be a loopless multigraph, let $P$ and $Q$ be edge-disjoint $a$--$b$ paths in $\Gamma$ of minimum total length, put $S=V(P)\cup V(Q)$, and direct both paths from $a$ to $b$. Then there is an integer function $\pi$ on $V(\Gamma)$ with $\pi(a)=0$ such that
\begin{enumerate}[label=\textup{(\roman*)}]
\item every arc $uv$ of $P$ or $Q$ satisfies $\pi(v)-\pi(u)\ge 1$, and
\item every edge $uv$ of $\Gamma$ which belongs to neither $P$ nor $Q$ satisfies $|\pi(v)-\pi(u)|\le 1$.
\end{enumerate}
Moreover, the vertices common to $P$ and $Q$ occur in the same order on both paths, and every integer is the label of at most two vertex occurrences on $P$ and $Q$, where a common vertex is counted twice.
\end{lemma}

\begin{proof}
Replace each physical edge by two opposite arcs, each of capacity one and cost one. The pair $P,Q$ is then a minimum-cost flow of value two from $a$ to $b$. Indeed, a cheaper flow could be chosen integral, and after cancelling opposite arcs on the same physical edge and deleting directed cycles it would give two $a$--$b$ paths using distinct physical edges, of smaller total length.

Since all costs are integers, the optimality conditions for minimum-cost flows~\cite[Chapter~9]{AMO93} supply integer potentials $\rho$ such that $\rho(v)\le\rho(u)+c$ for every arc $uv$ of cost $c$ in the residual network. A physical edge $uv$ which is not used by the flow yields both residual arcs $uv$ and $vu$ of cost one, so $|\rho(v)-\rho(u)|\le 1$. A used arc $uv$ yields the reverse residual arc $vu$ of cost $-1$, so $\rho(v)-\rho(u)\ge 1$. Hence $\pi=\rho-\rho(a)$ satisfies (i) and (ii).

The used arcs contain no directed cycle, since deleting such a cycle would reduce the cost without changing the value of the flow. If two common vertices occurred in opposite orders on $P$ and $Q$, the segments of $P$ and $Q$ between them would form a directed closed walk of used arcs, and hence would contain a directed cycle. Finally, labels strictly increase along each path by (i), so each integer is the label of at most one occurrence on each path. This completes the proof.
\end{proof}

We keep the notation of Lemma~\ref{lem:potential} and call $\pi(v)$ the \emph{label} of $v$ and $S$ the \emph{support} of the pair. A vertex of $S$ is \emph{common} if it lies on both paths and \emph{private} otherwise. An edge is \emph{used} if it belongs to $P\cup Q$ and \emph{unused} otherwise. Consecutive common vertices bound two internally disjoint subpaths, called the \emph{branches} between them, and their union is a cycle, called a \emph{support cycle}; two parallel edges form a support cycle of length two. These terms always refer to the whole pair, even when vertices are later omitted from a count.

Suppose that $\Gamma$ is simple. Then an unused edge with both ends in $S$ joins two private vertices on opposite branches of one support cycle. Indeed, if both ends of an unused edge lie on one of the paths, then their labels differ by at most one by (ii), so they are consecutive on that path by (i), and simplicity shows that the edge between them is the used one. In particular, a common vertex has no unused neighbour in $S$. Private vertices in different support cycles have a common vertex strictly between their labels, so their labels differ by at least two and they are not joined by an unused edge.

Whenever a pair of edge-disjoint paths is chosen below, we first minimise its total length, then maximise the sum of the degrees over its vertex occurrences, and finally maximise $|S|$. The degrees are measured in the fixed ambient simple graph before any core is contracted, and a contracted core receives a fixed score. We call a pair chosen in this way \emph{optimal}.

\begin{lemma}\label{lem:exchange}
Let $P,Q$ be an optimal pair with support $S$. Suppose that $a',v,b'$ are consecutive vertices of $P$, that $z\notin V(P)$, and that $a'z$ and $zb'$ are unused physical edges.
\begin{enumerate}[label=\textup{(\roman*)}]
\item If $v$ is private, then $\deg(z)\le\deg(v)$.
\item If $v$ is common and $z\notin S$, then $\deg(z)<\deg(v)$.
\end{enumerate}
The same holds with the roles of $P$ and $Q$ exchanged.
\end{lemma}

\begin{proof}
Replacing $a',v,b'$ by $a',z,b'$ on $P$ gives a path $P'$, since $z\notin V(P)$, and $P'$ is edge-disjoint from $Q$, since $a'z$ and $zb'$ are unused. The total length does not change, and the sum of the degrees over the vertex occurrences changes by $\deg(z)-\deg(v)$. If $v$ is private, then optimality gives $\deg(z)-\deg(v)\le 0$; this also covers the case that $z$ is a private vertex of $Q$, which becomes common. If $v$ is common and $z\notin S$, then a positive change contradicts the second priority. If $\deg(z)=\deg(v)$, then $v$ still occurs on $Q$ and $z$ is a new support vertex, so $|S|$ increases, contradicting the third priority. If $a'$ is a contracted core, then $a'z$ denotes a specified unused physical edge between the core and $z$; its end in the core need not be the end of $a'v$, but the replacement is still valid in the multigraph. This completes the proof.
\end{proof}

The next lemma separates neighbourhoods while discarding only a bounded number of support vertices. Put
\[
\gamma_0=1,\qquad\gamma_1=6,\qquad\gamma_2=41,\qquad\gamma_3=280.
\]
For an integer $j$, let $\partial_j$ be the set of ends of the used edges $uv$ with $\pi(u)\le j<\pi(v)$.

\begin{lemma}\label{lem:balls}
Let $\Gamma$, $P$, $Q$, $S$ and $\pi$ be as in Lemma~\ref{lem:potential}.
\begin{enumerate}[label=\textup{(\roman*)}]
\item For every integer $0\le t\le3$, every ball of radius $t$ in $\Gamma$ around a vertex of $S$ contains at most $\gamma_t$ vertices of $S$.
\item For every integer $j$ we have $|\partial_j|\le 4$ and $|S\cap N_\Gamma^2[\partial_j]|\le 164$. Moreover, if $x,y\in S\setminus N_\Gamma^2[\partial_j]$ satisfy $\pi(x)\le j<\pi(y)$, then $\dist_\Gamma(x,y)\ge 6$.
\item Suppose that $\Gamma$ is obtained from a simple graph $H$ by contracting a connected set $K$ to the vertex $a$. Then $|S\cap N_\Gamma^3[a]|\le 280$, and $N_H^2[Y]\cap N_H[K]=\emptyset$ for every $Y\subseteq S\setminus N_\Gamma^3[a]$.
\end{enumerate}
\end{lemma}

\begin{proof}
For (i), we first bound degrees inside the support. Let $v\in S$. If $v$ is private, then it has two used neighbours, and its unused neighbours in $S$ lie on the opposite path at the labels $\pi(v)-1$, $\pi(v)$ and $\pi(v)+1$: an unused edge to a vertex of the same path joins consecutive vertices by Lemma~\ref{lem:potential} and is parallel to a used edge. If $v$ is common, then it has at most four used neighbours, and for the same reason every unused edge from $v$ to $S$ is parallel to a used edge. Unused edges parallel to used edges create no new neighbours. Hence the underlying simple graph of $\Gamma[S]$ has maximum degree at most five.

Now let $x\in S$, let $y\in S$ with $\dist_\Gamma(x,y)\le t$, and fix a shortest $x$--$y$ path $R$. Either $y=x$, or the second vertex of $R$ is one of the at most five neighbours of $x$ in $S$, or $R$ first returns to $S$ after $\ell\ge 2$ edges. In the last case all these $\ell$ edges are unused, since used edges join two vertices of $S$, so by Lemma~\ref{lem:potential}(ii) the label of the first return differs from $\pi(x)$ by at most $\ell$. As every label carries at most two occurrences, the first return is one of at most $2(2\ell+1)=4\ell+2$ vertices of $S$, and the rest of $R$ has length at most $t-\ell$. By induction on $t$, the number of possible $y$ is at most
\[
1+5\gamma_{t-1}+\sum_{\ell=2}^{t}(4\ell+2)\gamma_{t-\ell}=\gamma_t\qquad(1\le t\le3),
\]
which gives the four displayed bounds directly.

For (ii), labels strictly increase along each of $P$ and $Q$, so each of them contains exactly one used edge $uv$ with $\pi(u)\le j<\pi(v)$ if $0\le j<\pi(b)$, and none otherwise. Hence $|\partial_j|\le 4$, and $|S\cap N_\Gamma^2[\partial_j]|\le 4\gamma_2=164$ by (i). Now let $x,y$ be as in (ii), and suppose that $\Gamma$ contains an $x$--$y$ path $R$ of length $L\le 5$. Both $x$ and $y$ have distance at least three from $\partial_j$, and a vertex of $R$ has distance at most two from $x$ or from $y$, so $R$ avoids $\partial_j$. Since $\pi(x)\le j<\pi(y)$, there are two consecutive visits $u,v$ of $R$ to $S$ with $\pi(u)\le j<\pi(v)$. Let $\ell$ be the length of the segment of $R$ between $u$ and $v$. Its internal vertices lie outside $S$, and if $\ell=1$ then the edge $uv$ is unused, because a used edge $uv$ would put $u,v$ into $\partial_j$. So all edges of this segment are unused and $\pi(v)-\pi(u)\le\ell$. Following a path among $P,Q$ through $u$ forwards, the labels increase by at least one at every step, so after at most $j-\pi(u)$ steps we reach the last vertex of label at most $j$, which lies in $\partial_j$. Similarly, following a path through $v$ backwards reaches a vertex of $\partial_j$ after at most $\pi(v)-j-1$ steps. Together with the parts of $R$ from $x$ to $u$ and from $v$ to $y$, these two walks show that
\[
\dist_\Gamma(x,\partial_j)+\dist_\Gamma(y,\partial_j)\le (L-\ell)+(\pi(v)-\pi(u)-1)\le L-1\le 4,
\]
contradicting $\dist_\Gamma(x,\partial_j)+\dist_\Gamma(y,\partial_j)\ge 6$.

For (iii), the first assertion follows from (i), since $a\in S$ and $\gamma_3=280$. Contraction preserves the distance between $K$ and every vertex outside $K$, so every vertex of $S\setminus N_\Gamma^3[a]$ has distance at least four from $K$ in $H$. A vertex within distance two of such a vertex therefore has distance at least two from $K$, and it does not lie in $N_H[K]$. This completes the proof.
\end{proof}

The separation in Lemma~\ref{lem:balls}(ii) concerns only two vertices on opposite sides of one cut. Hence it applies to several cuts simultaneously: if $Y_1,Y_2\subseteq S$ both avoid $N_\Gamma^2[\partial_j]$ and lie on opposite sides of $j$, then $\dist_\Gamma(Y_1,Y_2)\ge 6$, and therefore $N_\Gamma^2[Y_1]\cap N_\Gamma^2[Y_2]=\emptyset$. If $\Gamma$ is obtained from a simple graph $H$ by contraction, distances in $H$ are at least the corresponding distances in $\Gamma$, so we also have $N_H^2[Y_1]\cap N_H^2[Y_2]=\emptyset$.

\section{Removing long exact strips}\label{sec:strips}

We now shorten certain structures whose internal vertices have degree exactly $\delta$. These operations are carried out before any pair of paths is chosen and before any direction is prescribed; in particular, no direction is ever prescribed inside a part that is later removed. Throughout this section, $\delta\ge 3$ and $d=\delta+1$.

Let $m\ge 3$. An \emph{exact strip of length $m$} in a graph $H$ is a sequence $B_0,\dots,B_{m-1}$ of pairwise disjoint nonempty cliques of $H$, called bags, such that consecutive bags are completely joined and, for $1\le i\le m-2$ and every $z\in B_i$,
\[
N_H[z]=B_{i-1}\cup B_i\cup B_{i+1}\qquad\text{and}\qquad |B_{i-1}|+|B_i|+|B_{i+1}|=d .
\]

\begin{lemma}\label{lem:kernel}
Let $\delta\ge 3$, and let $G$ be a bridgeless graph of minimum degree at least $\delta$. Then there is a bridgeless graph $H$ of minimum degree at least $\delta$ with $|V(H)|\le|V(G)|$ which contains no exact strip of length $21$, such that every strong orientation $D_H$ of $H$ yields a strong orientation $D_G$ of $G$ with
\[
\diam(D_G)\le\diam(D_H)+\frac{3(|V(G)|-|V(H)|)}{d}+26 .
\]
\end{lemma}

The proof selects all replacements in the true-twin quotient before shortening any tube. We first give two local recovery constructions, each of which works for an arbitrary strong orientation of the smaller graph; the orientation induced on an individual replaced part need not be strong.

Let $\delta\ge 4$. A \emph{Type A cell} consists of a clique $C$ of order $\delta$, partitioned as $C=C_1\cup C_2$ with $|C_1|=p\ge 2$ and $|C_2|=q\ge 2$, together with two distinct vertices $l,r\notin C$, called its \emph{ports}, such that $l$ is the only neighbour outside $C$ of every vertex of $C_1$, and $r$ is the only neighbour outside $C$ of every vertex of $C_2$. Thus every vertex of $C$ has degree $\delta$. A \emph{Type A tube} of length $k$ is obtained by chaining $k$ cells, identifying the right port of each cell with the left port of the next one; each of the $k-1$ intermediate ports has degree $\delta$ and no further neighbours. The two end ports are \emph{retained}; they may have further neighbours, and they may coincide. The cliques of the cells and the intermediate ports form the \emph{interior} of the tube. If the $i$-th cell has the split $(p_i,q_i)$, then $q_i+p_{i+1}=\delta=p_i+q_i$, so $p_{i+1}=p_i$ and the split is the same along the tube.

Let $\delta\ge 3$. A \emph{Type B block} is an induced copy of $K_d$ minus one edge $lr$, all of whose vertices have degree $\delta$ in the whole graph. Its $\delta-1$ \emph{central} vertices have no neighbour outside the block, while each of its ports $l,r$ has exactly one. A \emph{Type B tube} of length $k$ is a sequence of $k$ blocks in which the right port of each block is joined to the left port of the next by its outside edge, called a \emph{joining edge}. The two end ports of the tube are retained and distinct, and they are the only vertices of the tube with edges leaving the tube; all other vertices of the tube form its interior. For a circular chain of blocks we retain one joining edge and cut the chain there.

A family of tubes is \emph{compatible} if the interiors are pairwise disjoint and no interior contains a retained vertex of another tube; retained vertices and edges may be shared. We replace a Type A tube of length $k\ge 4$ by a Type A tube of length two with the same split, and a Type B tube of length $k\ge 3$ by a Type B tube of length two. All new interior vertices are new, and every retained vertex and every edge between retained vertices is kept. The interiors have orders $dk-1$ and $dk-2$, respectively, so either replacement removes exactly $d(k-2)$ vertices. For a tube $T$ of length $k_T$ we put $\tau_T=3(k_T-2)$.

\begin{figure}[ht]
\centering
\begin{tikzpicture}[>=Stealth, vx/.style={circle,fill,inner sep=1.4pt}]
\node[vx,label=left:$l$] (l) at (0,0) {};
\node[vx,label=above:$u$] (u) at (1.3,1.0) {};
\node[vx,label=below:$v$] (v) at (1.3,-1.0) {};
\node[vx,label=above:$w$] (w) at (3.3,1.0) {};
\node[vx,label=below:$z$] (z) at (3.3,-1.0) {};
\node[vx,label=right:$x$] (x) at (2.3,0) {};
\node[vx,label=right:$r$] (r) at (4.6,0) {};
\draw[->] (l) -- (u);
\draw[->] (u) -- (w);
\draw[->] (w) -- (r);
\draw[->] (r) -- (z);
\draw[->] (z) -- (v);
\draw[->] (v) -- (l);
\draw[->] (u) -- (x);
\draw[->] (x) -- (v);
\draw[dashed,rounded corners=8pt] (0.75,-1.7) rectangle (3.85,1.7);
\node at (4.15,1.55) {$C$};
\node at (2.3,-2.3) {Type A cell};
\begin{scope}[xshift=7.6cm]
\node[vx,label=left:$l$] (L) at (0,0) {};
\node[vx,label=above:$h_1$] (H1) at (1.6,1.2) {};
\node[vx,label=below:$h_2$] (H2) at (1.6,-1.2) {};
\node[vx,label={[label distance=1pt]below right:$z$}] (Z) at (2.4,0) {};
\node[vx,label=right:$r$] (R) at (3.6,0) {};
\draw[->] (L) -- (H1);
\draw[->] (H1) -- (R);
\draw[->] (R) -- (H2);
\draw[->] (H2) -- (L);
\draw[->] (H1) -- (H2);
\draw[->] (L) -- (Z);
\draw[->] (Z) -- (R);
\draw[->] (H1) -- (Z);
\draw[->] (Z) -- (H2);
\node at (1.8,-2.3) {Type B block};
\end{scope}
\end{tikzpicture}
\caption{The prescribed orientation of a Type A cell (left), where $u,v\in C_1$, $w,z\in C_2$ and $x$ is any further vertex of $C_1\cup C_2$, and of a Type B block (right), where $z$ is any further central vertex. Edges which are not drawn are oriented arbitrarily.}
\label{fig:cells}
\end{figure}
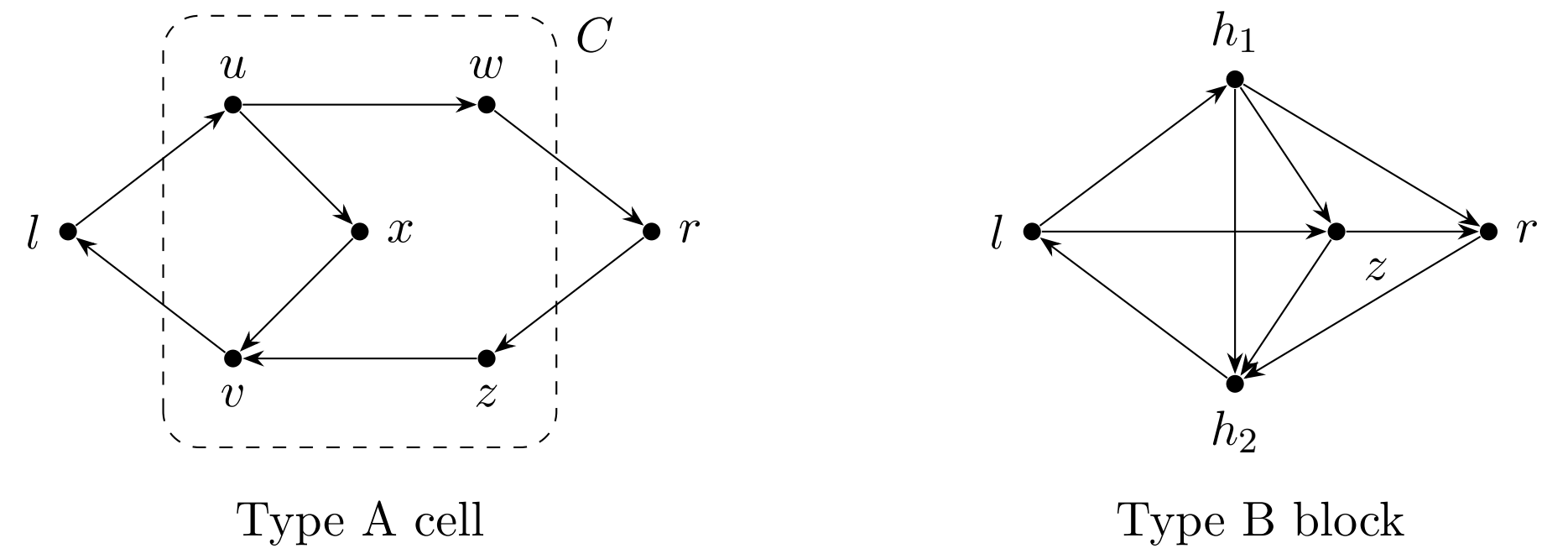

\begin{lemma}\label{lem:twoport}
Let $\delta\ge 3$, where Type A tubes are allowed only if $\delta\ge 4$. Let $G'$ be obtained from $G$ by the replacements in a compatible family $\mathcal{T}$ of Type A and Type B tubes. Then every strong orientation $D'$ of $G'$ yields a strong orientation $D$ of $G$ which agrees with $D'$ on every retained edge and satisfies
\[
\diam(D)\le\diam(D')+\sum_{T\in\mathcal{T}}\tau_T+12=\diam(D')+\frac{3(|V(G)|-|V(G')|)}{d}+12 .
\]
Moreover, the replacements preserve simplicity, connectivity, bridgelessness, the minimum degree being at least $\delta$, and the degree of every retained vertex.
\end{lemma}

\begin{proof}
We first fix the orientation inside the expanded tubes; this is done before any ordered pair is considered.

Let $T$ be a Type A tube of length $k\ge 4$. In each cell with ports $l,r$, choose distinct $u,v\in C_1$ and $w,z\in C_2$ and prescribe
\[
l\to u\to w\to r\to z\to v\to l ,
\]
together with $u\to x\to v$ for every other vertex $x$ of $C_1\cup C_2$; the remaining edges of the cell are oriented arbitrarily (see Figure~\ref{fig:cells}). These prescriptions do not conflict. In one cell, each port reaches every vertex of $C_1\cup C_2$, and is reached from it, within five steps; both port-to-port distances are three; and any two vertices of $C_1\cup C_2$ are at distance at most six, since the six-cycle above handles its own vertices and every other vertex is entered from $u$ and leaves to $v$ in one step. Give the left end of the tube the coordinate $0$, the $j$-th intermediate port the coordinate $3j$, the right end the coordinate $3k$, and every vertex of the clique of the $j$-th cell the coordinate $3(j-1)+3/2$. Both crossings of the tube have length $3k=\tau_T+6$. A vertex of coordinate $t$ has directed distance at most $t+4$ to and from the left end and at most $3k-t+4$ to and from the right end; for intermediate ports the additive $4$ can be omitted. Any two interior vertices are joined by a directed path of length at most $3k+4=\tau_T+10$: use at most five steps in each of the two end cells and three steps in every cell in between, or the bound six inside one cell. For each interior vertex we choose, once and for all, a nearer end. Its distances to and from this end are at most $3k/2+4\le\tau_T+4$, since $k\ge 4$, and its distances to and from either end are at most $\tau_T+10$. These statements remain true when the two ends coincide.

Let $T$ be a Type B tube of length $k\ge 3$. In each block with ports $l,r$, choose central vertices $h_1,h_2$, let $Z$ be the set of the remaining central vertices, and prescribe
\[
l\to h_1\to r\to h_2\to l,\qquad h_1\to h_2,\qquad l\to z\to r\quad\text{and}\quad h_1\to z\to h_2\quad(z\in Z),
\]
orienting the edges inside $Z$ arbitrarily. Both port-to-port distances in a block are two, and the block has diameter at most three; for instance, $z\to h_2\to l\to z'$ for $z,z'\in Z$, and the four vertices $l,h_1,h_2,r$ give the same bound for all other types of pairs.

The joining edges of $T$ are oriented according to $D'$ as follows. In the compressed tube in $G'$, every proper initial segment of blocks has exactly two edges to the rest of $G'$, namely the outside edge at its first port and its last joining edge. Since $D'$ is strong, these two edges point in opposite directions across this cut. Applying this to successive initial segments shows that all joining edges of the compressed tube point in the same direction along the tube, and we order the tube so that this direction is forwards. Then the outside edge at the first end port points into the tube. The same cut argument, applied to the last block, shows that the outside edge at the last end port points out of the tube; for a circular chain, the retained joining edge plays both roles. In particular, $D'$ contains no directed path crossing the compressed tube backwards. Note that this uses the strong connectivity of the whole of $D'$, not of a single block. We orient all joining edges of the expanded tube forwards and use the block orientations fixed above. Let $a$ and $b$ be the first and last end ports. The forward crossing has length $3k-1=\tau_T+5$. If $x$ lies in the $i$-th block and $y$ in the $j$-th block, then
\[
\dist_D(x,b)\le 3(k-i)+3\qquad\text{and}\qquad\dist_D(a,y)\le 3(j-1)+3 .
\]
In particular $\dist_D(x,b)\le\tau_T+6$ and $\dist_D(a,x)\le\tau_T+6$ for every interior vertex $x$. If $i>j$, the sum of the two displayed bounds is at most $\tau_T+6$. If $i<j$, a forward route from $x$ to $y$ has length at most $3(j-i)+4\le\tau_T+7$, and if $i=j$ the block bound three applies.

All retained edges keep their directions from $D'$. This completes the definition of $D$.

We next lift directed paths of $D'$, charging each tube at most once. A shortest directed path of $D'$ between retained vertices meets a two-port tube in at most one crossing, since returning to the entry port, or crossing a second time, would repeat a port. A tube whose two end ports coincide is never entered by such a path. A crossing of a compressed Type A tube has at least six edges, and a crossing of a compressed Type B tube has at least five. Replacing a crossing by the corresponding expanded crossing therefore increases the length by at most $\tau_T$.

Now let $x,y$ be vertices of $G$ which do not lie in the interior of a common tube. We map an interior vertex of a Type A tube to its chosen nearer end, a vertex $x$ in the interior of a Type B tube to the last end port of that tube, and a vertex $y$ in the interior of a Type B tube to the first end port; retained vertices are mapped to themselves. Let $R$ be a shortest directed path of $D'$ from the image of $x$ to the image of $y$. Suppose that $x$ lies in a Type A tube $T$ and that $R$ crosses $T$. This crossing is the initial segment of $R$, since a later crossing would revisit the first vertex of $R$. We delete it, which removes at least six edges, and access the other end of $T$ directly from $x$, at cost at most $\tau_T+10$ instead of $\tau_T+4$. So in all cases the tube containing $x$ costs at most $\tau_T+4$. The vertex $y$ is treated symmetrically, and if both $x$ and $y$ lie in Type A tubes, the two deleted crossings are disjoint because the interiors are different. A Type B tube containing $x$ or $y$ is not met by $R$ at all: the image of $x$ is the last port of its tube, the image of $y$ is the first port of its tube, a backward crossing does not exist, and a forward crossing would revisit the first or the last vertex of $R$. Consequently, each tube containing $x$ or $y$ contributes its $\tau_T$ once and at most four further steps for Type A or six for Type B, and every other tube met by $R$ contributes at most its $\tau_T$. This proves the bound on $\dist_D(x,y)$ for such pairs.

If $x$ and $y$ lie in the interior of the same Type A tube, the internal route has length at most $\tau_T+10$. If they lie in the interior of the same Type B tube, the forward order and the case of a single block were treated above. In the backward order, that is, when $x$ lies in a later block than $y$, we go from $x$ to the last end port $b$, follow a shortest directed path of $D'$ from $b$ to the first end port $a$, and go from $a$ to $y$. The path from $b$ to $a$ does not meet the tube, since a backward crossing does not exist and a forward crossing would revisit $a$ or $b$. The two accesses cost at most $\tau_T+6$ together, and the other tubes met by the path are lifted as before. Since every ordered pair is joined by a directed path, $D$ is strong.

Finally, we check the properties of $G'$. The sizes of the groups adjacent to a retained Type A port are unchanged, and a retained Type B port still has $\delta-1$ central neighbours and the same outside edge. Hence the degrees of retained vertices are unchanged, and every new vertex has degree $\delta$. The new interiors consist of new vertices, so $G'$ is simple. Each replacement joins exactly the same retained vertices as before, which preserves connectivity, and also connectivity after the deletion of any retained edge. Every new edge of a Type A tube and every edge inside a new Type B block lies in a triangle. For a Type B tube, deleting an original joining edge and using that $G$ is bridgeless gives a path outside the tube between its two end ports, because such a path has to leave the tube through one end and return through the other. Replacing other tubes preserves such an outside path, and it closes the new chain of blocks into a cycle through every new joining edge. Thus $G'$ is bridgeless. This completes the proof.
\end{proof}

The remaining replacements use clique boundaries. A \emph{clique-layer tube} has bags $V_0,\dots,V_L$, each a clique of order at least two. Consecutive bags are completely joined. The interior bags $V_1,\dots,V_{L-1}$ are pairwise disjoint and disjoint from the boundary bags $V_0,V_L$. Every vertex $z$ of an interior bag $V_j$ satisfies $N[z]=V_{j-1}\cup V_j\cup V_{j+1}$, where
\[
|V_{j-1}|+|V_j|+|V_{j+1}|=d\qquad(1\le j\le L-1).
\]
The two boundary bags may coincide, and they may have arbitrary further edges. Taking differences of consecutive equations shows that the bag sizes are periodic with period three. Choose $L'\in\{6,7,8\}$ with $L'\equiv L\pmod 3$, and assume that $L\ge 2L'$. We replace the interior by $L'-1$ new interior bags with the same periodic sizes, keeping both boundary bags, so that the new tube has $L'$ interfaces between consecutive bags. For such a tube $T$ put $\tau_T=L-L'$. The number of removed vertices is exactly $\tau_Td/3$, and the sizes of the bags next to each boundary bag are unchanged. Compatibility again means that the interiors are disjoint and that no boundary vertex lies in another interior.

\begin{lemma}\label{lem:layer}
Let $G'$ be obtained from $G$ by the replacements in a compatible family $\mathcal{T}$ of clique-layer tubes. Then every strong orientation $D'$ of $G'$ yields a strong orientation $D$ of $G$ which agrees with $D'$ on every retained edge and satisfies
\[
\diam(D)\le\diam(D')+\sum_{T\in\mathcal{T}}\tau_T+14=\diam(D')+\frac{3(|V(G)|-|V(G')|)}{d}+14 .
\]
Moreover, the replacements preserve simplicity, connectivity, bridgelessness, the minimum degree being at least $\delta$, and the degree of every retained vertex.
\end{lemma}

\begin{proof}
We first fix the expanded orientation and its local routes. In every interior bag $V_j$, choose distinct vertices $u_j,v_j$ and put $Z_j=V_j\setminus\{u_j,v_j\}$. Prescribe
\[
u_j\to v_j\quad\text{and}\quad u_j\to z\to v_j\quad(z\in Z_j)\qquad\text{for }1\le j\le L-1,
\]
\[
u_j\to u_{j+1},\qquad v_{j+1}\to v_j\qquad\text{and}\qquad v_j\to u_{j+1}\qquad\text{for }1\le j\le L-2,
\]
and, for every $a\in V_0$ and every $b\in V_L$,
\[
a\to u_1,\qquad v_1\to a,\qquad u_{L-1}\to b\qquad\text{and}\qquad b\to v_{L-1}.
\]
Keep the directions of all retained edges and orient the remaining edges arbitrarily. These prescriptions remain consistent when the boundary bags coincide, since the first and the last interior bags are distinct.

The $u$-lane and the $v$-lane give both boundary crossings, from any $a\in V_0$ to any $b\in V_L$ and back, in $L$ steps. Two vertices $a,a'$ of the left boundary bag are joined by the route $a\to u_1\to v_1\to a'$ of length three, and two vertices $b,b'$ of the right boundary bag by the route $b\to v_{L-1}\to v_{L-2}\to u_{L-1}\to b'$ of length four. For every $x\in V_j$ with $1\le j\le L-1$ and all $a\in V_0$ and $b\in V_L$, the same arcs give
\[
\dist_D(a,x),\ \dist_D(x,a)\le j+1\qquad\text{and}\qquad\dist_D(b,x),\ \dist_D(x,b)\le L-j+3 .
\]
For the right-hand estimates, note that every vertex of $V_j$ reaches $u_j$, and is reached from $v_j$, within three steps: the turn from $v_j$ to $u_j$ uses $v_j\to v_{j-1}\to u_j$ if $j\ge 2$ and $v_1\to a\to u_1$ if $j=1$, and a vertex of $Z_j$ uses its arc to $v_j$ or from $u_j$. The same observations give, for vertices of the bags $V_j$ and $V_{j'}$ with $1\le j,j'\le L-1$, a route of length at most $|j-j'|+4\le L+2$. For every interior vertex we choose, once and for all, a nearer boundary side and a vertex of the boundary bag on that side. Both accesses between the interior vertex and this boundary vertex have length at most $L/2+3\le\tau_T+3$, since $L\ge 2L'$. The access between an interior vertex and any vertex of either boundary bag has length at most $L+3=\tau_T+L'+3$. Note also that $L\ge 2L'$ gives $\tau_T\ge L'\ge 6$.

Next we deal with repeated visits of a path to a tube. A \emph{visit} of a path of $D'$ to a tube is a maximal subpath whose internal vertices lie in the interior of the tube and whose ends are retained; its entry side and its exit side are the sides of the first and of the last interior bag that it meets. A visit is a \emph{crossing} if these sides are different, and a \emph{return} otherwise. A crossing has at least $L'$ edges, and a return has at least two. These definitions use the interior bags, and therefore they also apply when the boundary bags coincide.

Let $R$ be a directed path of $D'$ between retained vertices. We process $R$ from its beginning and replace disjoint intervals of $R$. Suppose that we have reached the first remaining visit to a tube $T$. If the rest of $R$ contains a crossing of $T$, we replace the interval of $R$ from the beginning of the current visit to the end of the last visit of $R$ to $T$ by a walk inside the expanded tube. If the two ends of this interval lie on opposite sides, the new walk has length at most $L$, whereas the interval contains a crossing and hence at least $L'$ edges, so the length increases by at most $\tau_T$. If the two ends lie on the same side, a walk of length at most four replaces an interval with at least $L'\ge 6$ edges, and the length does not increase. In both cases no visit to $T$ remains. If the rest of $R$ contains no crossing of $T$, then all remaining visits to $T$ are returns. We replace the interval from the beginning of the current visit to the end of the last visit on the same side by a walk of length at most four on that side, which increases the length by at most $4-2=2$. Afterwards only returns on the other side can remain, and one further replacement of the same kind removes them at a cost of at most two. In this case the total increase is at most $4\le\tau_T$. Note that visits to the two sides are not automatically replaced by a crossing, since the part of $R$ between them may use a short route outside the tube; the procedure above makes no assumption about distances outside the tube. Each tube is charged at most $\tau_T$, once, and deleting an interval may remove visits to other tubes but never creates a new visit.

Finally, consider vertices $x,y$ of $G$ which do not lie in the interior of a common tube. We map an interior vertex to its chosen boundary vertex and a retained vertex to itself, and we take a shortest directed path $R$ of $D'$ between the images. Suppose that $x$ lies in the interior of a tube $T$. If $R$ contains a crossing of $T$, we delete the initial part of $R$ up to the end of the last visit of $R$ to $T$ and go from $x$ directly to the end of that visit. At least $L'$ edges disappear and the new access costs at most $\tau_T+L'+3$, so the net cost is at most $\tau_T+3$, and no visit to $T$ remains. If $R$ contains no crossing of $T$, we use the access to the chosen boundary vertex, of length at most $\tau_T+3$, and leave the remaining returns to the procedure above. The vertex $y$ is treated symmetrically on the remaining path: if a crossing of its tube remains, we delete the final part of $R$ from the beginning of the first visit to that tube and access $y$ directly. The deleted initial and final parts are disjoint, also when boundary vertices are shared, and the first deletion may remove a crossing of the other tube. A tube containing $x$ or $y$ whose visits remain has only returns, which cost at most four more. Thus each tube containing $x$ or $y$ contributes its $\tau_T$ and at most seven further steps, and all other tubes are charged as before. This proves the bound. If $x$ and $y$ lie in the interior of the same tube, we use the internal route of length at most $L+2\le\tau_T+10$. All these walks use the orientation fixed at the beginning, so $D$ is strong.

The periodic sizes give the exact number of removed vertices and unchanged degrees of retained vertices, and every new interior vertex has degree $\delta$. The new interiors consist of new vertices, so the graph stays simple. Every new edge incident with an interior vertex lies in a triangle, because every bag has order at least two. Each replacement joins the same boundary vertices as before, which preserves connectivity, and also connectivity after the deletion of any retained edge. Thus connectivity and bridgelessness are preserved. This completes the proof.
\end{proof}

\begin{proof}[Proof of Lemma~\ref{lem:kernel}]
Two vertices are \emph{true twins} if they have the same closed neighbourhood. Their equivalence classes partition $V(G)$ into cliques, and two distinct classes are either completely joined or have no edge between them. Let $\widehat G$ be the simple quotient graph, with each class weighted by its order. A class is \emph{marked} if it has degree two in $\widehat G$ and its vertices have degree exactly $\delta$ in $G$; singleton classes are allowed. The components induced by the marked classes are paths and cycles. A cycle component is the whole of $\widehat G$, since each of its classes already has both quotient neighbours on the cycle. It cannot be a triangle, whose three classes would have identical closed neighbourhoods.

Write $C_0,\ldots,C_{k-1}$ for the classes along a marked path. For $1\le i\le k-2$, the two quotient neighbours of $C_i$ are $C_{i-1}$ and $C_{i+1}$, and hence
\[
N_G[z]=C_{i-1}\cup C_i\cup C_{i+1},\qquad
|C_{i-1}|+|C_i|+|C_{i+1}|=d\quad(z\in C_i).
\]
The same equations hold cyclically on a marked cycle. Consecutive equations show that the class sizes have period three. We select one tube in each marked path with $k\ge17$ and each marked cycle with $k\ge16$, as follows.

\emph{All three periodic sizes are at least two.} On a path, keep its first and last classes as boundary bags and use the intervening clique-layer tube, with $L=k-1\ge16$ interfaces. On a cycle, keep one class and cut there, so that the boundary bags coincide and $L=k\ge16$. In either case, choose $L'\in\{6,7,8\}$ congruent to $L$ modulo three. Then $L\ge2L'$, so Lemma~\ref{lem:layer} applies.

\emph{Exactly one periodic size is one.} The other two sizes, say $p,q\ge2$, satisfy $p+q=\delta$. Between consecutive singleton classes, the two intervening classes form a Type A cell: their union is a clique of order $\delta$, and the two singleton classes are its ports. On a path, take the tube from the first singleton to the last singleton. At most two classes precede it and at most two follow it. If it has $h$ cells, it contains $3h+1$ classes, so $k\le3h+5$ and $k\ge17$ implies $h\ge4$. On a cycle, periodicity implies that $k$ is divisible by three. Cutting at a singleton class gives a Type A tube with coinciding end ports and $h=k/3\ge6$ cells. We shorten the chosen tube to two cells.

\emph{Exactly two periodic sizes are one.} The remaining size is $\delta-1$. A class of that size together with its two adjacent singleton classes is a Type B block whenever all three classes lie in the component. All its vertices have degree $\delta$, and its two ports are not adjacent. Successive such blocks are joined by their actual joining edges. On a path, take all these complete blocks as one tube. At most two classes precede this tube and at most two follow it; if it has $h$ blocks, then $k\le3h+4$, so $h\ge5$ when $k\ge17$. On a cycle, retain one joining edge and cut there, obtaining $h=k/3\ge6$ blocks. We shorten the chosen tube to two blocks.

These cases exhaust the possibilities, since three singleton sizes would give $d=3$, whereas $\delta\ge3$. Every selected tube uses only classes of its own marked component. Thus different selected tubes have disjoint vertex sets, and both the Type A/B family and the clique-layer family are compatible. We choose both families in $G$ before performing any replacement. The interior of a selected tube has no neighbours outside that tube, and the two families have disjoint boundary vertices as well. Consequently, replacing the Type A/B family leaves every selected clique-layer tube intact. Apply the replacements of Lemma~\ref{lem:twoport}, then those of Lemma~\ref{lem:layer}, and call the resulting graph $H$.

We verify the true-twin classes of $H$ directly. Each new bag is a clique of true twins. The quotient path inserted by a replacement has at least five interfaces; a circular replacement has at least six classes. Distinct adjacent new bags are distinguished by their outer neighbouring bags. A new end bag and its retained boundary class are distinguished by the other quotient neighbour of the boundary class; in a circular replacement this neighbour is the bag or port at the other end. Each retained boundary class has a fresh neighbouring bag that is adjacent to no other retained class, so it cannot merge with another retained class. Every other retained class has no neighbour in a deleted interior, and therefore its old distinctions from the other such classes persist. Nonadjacent vertices cannot be true twins, and no retained class splits, since all of its vertices receive the same new neighbours. The true-twin classes of $H$ are therefore exactly the retained classes and the new bags.

Every retained class keeps its weight, quotient degree and vertex degree, and hence its marked status. Every new bag is marked. The marked components can thus only be shortened, without merging. A shortened marked path has at most eleven classes in the Type A case, ten in the Type B case, and nine in the clique-layer case; the first two bounds include the at most two unchanged classes at each end. A shortened marked cycle has six classes in the Type A/B cases and at most eight in the clique-layer case. Untouched paths have at most sixteen classes and untouched cycles at most fifteen. Thus $H$ has no marked component containing seventeen classes.

Suppose that $B_0,\ldots,B_{20}$ is an exact strip in $H$. Each $B_i$ with $1\le i\le19$ is a full true-twin class. Indeed, its vertices have the same closed neighbourhood, and an additional twin would have to lie in an adjacent bag. The nonempty bags at distance two distinguish adjacent bags; at the two ends, the exact neighbourhoods of $B_2$ and $B_{18}$ give the same distinction. It follows that $B_2,\ldots,B_{18}$ are seventeen consecutive marked classes, a contradiction. Hence $H$ contains no exact strip of length $21$.

Both replacement lemmas preserve simplicity, connectivity, bridgelessness and the minimum-degree condition. Given a strong orientation of $H$, first recover through the clique-layer replacements and then through the Type A/B replacements. The order terms add up to $3(|V(G)|-|V(H)|)/d$, and the additive terms are at most $14+12=26$. Each recovery produces one orientation serving all ordered pairs. This proves the lemma.
\end{proof}

In all subsequent counts we work in the fixed reduced graph before any core is contracted. New vertices created by the reduction are not identified with deleted vertices of the original graph; their contribution has already been accounted for by the exact order difference in Lemma~\ref{lem:kernel}. Performing the reduction first also avoids assuming that it preserves an orientation prescribed inside a part that it removes.

\section{Strict counting and extension of an oriented core}\label{sec:counting}

Throughout this section, degrees and neighbourhoods refer to a fixed simple graph $G$, before any contraction. Unless weaker hypotheses are stated explicitly, $G$ has minimum degree at least $\delta\ge 3$, and $d=\delta+1$ is kept after contraction. Pairs of edge-disjoint paths are optimal pairs in the sense of Section~\ref{sec:paths}, and $\pi$ denotes the labels of Lemma~\ref{lem:potential}. Thus labels strictly increase along each path, an unused edge has label difference at most one, and every label has at most two occurrences. Recall that, in a simple graph, an unused edge inside the support joins two private vertices on opposite branches of one support cycle, and that a common vertex has no unused neighbour in the support.

Choose $s$ support cycles, called \emph{chosen} cycles, and designate one branch of each of them as \emph{positive}; the other branch is called the \emph{zero branch}. Assign multiplicity two to every common vertex and to every private vertex on a positive branch, multiplicity zero to every private vertex on a zero branch, and multiplicity one to all remaining private vertices. Write $w(v)$ for this multiplicity. For a set $Y$ of support vertices, whose elements we call \emph{retained} while the other support vertices are called \emph{discarded}, put
\[
T_Y=\sum_{v\in Y}w(v)\qquad\text{and}\qquad h=|N_G[Y]|,
\]
and let $b$ be the number of used edges with exactly one end in $Y$. Retained vertices of multiplicity zero remain members of $Y$.

Let $t_0<\dots<t_{\Lambda-1}$ be the labels of the vertices of $Y$, called the \emph{retained labels}. Every vertex $z$ of $G$ receives the \emph{compressed label} $\lambda(z)=\max\{i:t_i\le\pi(z)\}$, where we put $\lambda(z)=0$ if $\pi(z)<t_0$. An unused edge still has compressed span at most one. We now augment the multiplicities so that every retained label has total multiplicity two. A retained label with total multiplicity less than two is called \emph{deficient}. A deficient label outside the open label intervals of the chosen cycles has exactly one retained representative, which is private of multiplicity one, and we increase its multiplicity to two. A deficient label inside a chosen cycle has exactly one retained representative, which lies on the zero branch, and we increase its multiplicity from zero to two. These are the only possibilities. We write $a(v)$ for the augmented multiplicities, extended by zero outside $Y$, and put
\[
k=\sum_{v\in Y}\bigl(a(v)-w(v)\bigr)\qquad\text{and}\qquad U=\sum_{v\in Y}a(v)\bigl(\deg_G(v)-\delta\bigr).
\]
Then $T_Y+k=2\Lambda$.

Let $\Omega$ be the set of used edges with exactly one end in $Y$ together with the used edges inside $Y$ whose compressed span is at least two. The edges of $\Omega$ are omitted only in the following incidence count. For $z\in N_G[Y]$ put
\[
\ell(z)=a(z)+\sum_{zw\in E(G)\setminus\Omega}a(w)\qquad\text{and}\qquad D_0=\sum_{z\in N_G[Y]}\bigl(6-\ell(z)\bigr).
\]
We call the vertices of $N_G[Y]$ \emph{centres} and $\ell(z)$ the \emph{load} of $z$. All contributions to $\ell(z)$ come from three consecutive compressed labels, so $0\le\ell(z)\le 6$. Let $S$ be the full support, and define
\[
\Phi=\{z\in N_G[Y]\setminus S:\ \ell(z)=6\},\qquad J^*=\bigcup_{z\in\Phi}\bigl(N_G(z)\setminus N_G[Y]\bigr)\qquad\text{and}\qquad J=|J^*| .
\]
Every vertex of $\Phi$ has neighbours of positive multiplicity at three distinct labels. In particular $J^*\subseteq N_G^2[Y]\setminus N_G[Y]$. For $Y=\emptyset$ all these quantities are zero.

\begin{lemma}\label{lem:strict}
With the notation above, the following hold.
\begin{enumerate}[label=\textup{(\roman*)}]
\item If every vertex of $Y$ has degree at least three, then
\[
\sum_{v\in Y}w(v)\bigl(\deg_G(v)+1\bigr)\le 6h+4s+2b .
\]
\item Suppose that $G$ contains no exact strip of length $21$, and that either $\delta\ge4$ or $\delta=3$ and $k=0$. Put $\kappa=1470$. Then
\[
\frac{T_Y}{2}\le\frac{3(h+J/2)}{d}-\frac{|Y|}{2\kappa d}+\frac{1764s+588b+628}{2\kappa d}+\frac{2s+b}{d}.
\]
\end{enumerate}
Part~\textup{(i)} uses neither a degree bound outside $Y$ nor the second and third priorities in the choice of the pair.
\end{lemma}

\begin{proof}
We may assume that $Y\ne\emptyset$. Let $E_a$ be the augmented incidence omitted on $\Omega$: an edge $uv\in\Omega$ inside $Y$ contributes $a(u)+a(v)$, and an edge of $\Omega$ leaving $Y$ contributes the multiplicity of its end in $Y$. Let $K_0$ be the number of deficient labels outside the open label intervals of the chosen cycles, and for a chosen cycle with positive branch $A$ let $K_A$ be the number of its deficient labels. Then $k=K_0+2\sum_AK_A$, and we claim that
\[
E_a\le 4k+4s+2b .
\]
Consider an unchosen cycle. An edge of $\Omega$ inside $Y$ on one branch of it skips a retained label which is represented only on the opposite branch, and this label is deficient. Different edges skip different labels, since the open label intervals of the used edges on one branch are disjoint, and edges on opposite branches cannot both skip a retained label. Each such edge contributes at most four, so these edges contribute at most $4K_0$ in total. In a chosen cycle with positive branch $A$, an edge of $\Omega$ inside $Y$ on $A$ skips a retained label represented only on the zero branch, which is deficient, and as before these edges contribute at most $4K_A$. On the zero branch, the original multiplicities are positive only at the two common ends, which contribute at most four in total, and each augmented representative adds at most four through its two used edges. So the zero branch contributes at most $4K_A+4$. Finally, each edge of $\Omega$ leaving $Y$ contributes at most two. Summing gives $E_a\le 4K_0+\sum_A(8K_A+4)+2b\le 4k+4s+2b$.

Every vertex $v\in Y$ is counted in $\ell(z)$ for $z=v$ and for every neighbour $z$ of $v$ with $zv\notin\Omega$. Hence
\[
\sum_{v\in Y}a(v)\bigl(\deg_G(v)+1\bigr)=\sum_{z\in N_G[Y]}\ell(z)+E_a=6h-D_0+E_a .
\]
Put $W=\sum_{v\in Y}w(v)(\deg_G(v)+1)$. Since every augmented vertex has degree at least three, the left-hand side is at least $W+4k$. Therefore $W+4k\le 6h+E_a\le 6h+4k+4s+2b$, which proves (i).

For (ii), the left-hand side of the identity above equals $d(T_Y+k)+U$, since every retained label has total multiplicity two. Put
\[
\Psi=D_0+U+(d-4)k+3J ,
\]
a sum of nonnegative terms. Using $E_a\le 4k+4s+2b$, we obtain $dT_Y\le 6(h+J/2)-\Psi+4s+2b$, that is,
\[
\frac{T_Y}{2}\le\frac{3(h+J/2)}{d}-\frac{\Psi}{2d}+\frac{2s+b}{d}.
\]
It remains to bound the number $\Lambda$ of retained labels in terms of $\Psi$.

We mark the following compressed labels: the label of every centre $z\in N_G[Y]$ with $\ell(z)<6$; the label of every vertex of $J^*$; the label of every retained vertex $v$ with $a(v)>0$ and $\deg_G(v)>\delta$; every deficient label; the label of every end $v$ with $a(v)>0$ of an edge of $\Omega$; the labels of the two common ends of every chosen cycle, whether these ends are retained or not; and the first and the last retained label. The number of marked labels is at most
\[
D_0+J+U+k+E_a+2s+2\le D_0+J+U+5k+6s+2b+2 .
\]
A retained label is \emph{bad} if it is within three compressed units of a marked label, and \emph{clean} otherwise. The number $\Nbad$ of bad labels satisfies
\[
\Nbad\le 7(D_0+J+U+5k+6s+2b+2).
\]

For every $i$ let $C_i=\{z\in N_G[Y]:\lambda(z)=i\}$. A vertex $z\in C_i$ is called \emph{exact} if
\[
N_G[z]=C_{i-1}\cup C_i\cup C_{i+1}\qquad\text{and}\qquad|C_{i-1}|+|C_i|+|C_{i+1}|=d .
\]
We now show that, for every clean label $j$, every vertex of $C_j$ is exact. Fix a clean label $j$. Then $4\le j\le\Lambda-5$, and no label $i$ with $|i-j|\le 3$ is marked. Hence, for every such $i$, every centre in $C_i$ has load six; the label $t_i$ is not deficient, so $a=w$ at $t_i$ and the retained vertices at $t_i$ have total original multiplicity two; no vertex of $J^*$ has compressed label $i$; every retained vertex of positive multiplicity at $t_i$ has degree $\delta$ and is incident with no edge of $\Omega$; and $t_i$ is not the label of a common end of a chosen cycle.

We use the following consequence of load six several times. Let $z$ be a centre with $\ell(z)=6$ and $|\lambda(z)-j|\le 3$. An edge which contributes positively to a load has an end in $Y$ of positive multiplicity and does not belong to $\Omega$. Such an edge is either unused, or a used edge inside $Y$ which is not in $\Omega$; in both cases its compressed span is at most one. Since every retained label has total multiplicity two, $z$ collects multiplicity only from the labels $\lambda(z)-1$, $\lambda(z)$ and $\lambda(z)+1$, and it must collect all of it. Hence $z$ is equal or adjacent to every retained vertex of positive multiplicity whose compressed label differs from $\lambda(z)$ by at most one.

First, let $v$ be a retained vertex of positive multiplicity with $\lambda(v)=i$, where $|i-j|\le 2$. Every vertex of $C_{i-1}\cup C_i\cup C_{i+1}$ is a centre of load six whose compressed label differs from $i$ by at most one, so it is equal or adjacent to $v$. Conversely, $v$ is incident with no edge of $\Omega$, so every neighbour of $v$ lies in $N_G[Y]$ and has compressed label in $\{i-1,i,i+1\}$. Hence $N_G[v]=C_{i-1}\cup C_i\cup C_{i+1}$, and since $\deg_G(v)=\delta$, this set has exactly $d$ vertices. As every retained label carries a vertex of positive multiplicity, we conclude that
\[
|C_{i-1}|+|C_i|+|C_{i+1}|=d\qquad\text{whenever }|i-j|\le 2,
\]
and that every retained vertex of positive multiplicity in such a set $C_i$ is exact. Note that this conclusion about the sizes of the sets $C_i$ is obtained before any degree claim about other vertices.

Second, let $z\in C_i\setminus S$ with $|i-j|\le 2$. Then $z$ is a centre of load six, so $z\in\Phi$. Since $a(z)=0$ and unused edges change labels by at most one, the neighbours of $z$ of positive multiplicity lie at the labels $\pi(z)-1$, $\pi(z)$ and $\pi(z)+1$, and they carry total multiplicity six. Hence these three integers are retained labels, and $\pi(z)=t_i$, $t_{i-1}=t_i-1$ and $t_{i+1}=t_i+1$. Let $v$ be a vertex of positive multiplicity at the label $t_i$. Suppose first that $v$ is private, and let $v^-$ and $v^+$ be its neighbours on its path. The edges $vv^-$ and $vv^+$ are not in $\Omega$, so $v^-$ and $v^+$ are retained and have compressed labels $i-1$ and $i+1$. They have positive multiplicity: this is clear in an unchosen cycle, and in a chosen cycle $v$ lies on the positive branch while $v^-,v^+$ are not common ends of the cycle, because these ends are marked. By the consequence of load six, $z$ is adjacent to $v^-$ and $v^+$, and these edges are unused since $z\notin S$. Replacing $v$ by $z$ on its path is the exchange of Lemma~\ref{lem:exchange}(i), which gives $\deg_G(z)\le\deg_G(v)=\delta$; hence $\deg_G(z)=\delta$. Suppose next that $v$ is common. In the same way, the neighbours of $v$ on $P$ are retained vertices of positive multiplicity at the labels $t_i-1$ and $t_i+1$, and $z$ is adjacent to both of them. Lemma~\ref{lem:exchange}(ii) then gives $\deg_G(z)<\deg_G(v)=\delta$, which is impossible. So the middle label of $z$ carries no common vertex, and $\deg_G(z)=\delta$. Every edge at $z$ is unused and thus not in $\Omega$. A neighbour of $z$ outside $N_G[Y]$ would belong to $J^*$ and would have compressed label within one of $i$, which is not marked. Hence all neighbours of $z$ lie in $C_{i-1}\cup C_i\cup C_{i+1}$. This set has $d$ vertices and $\deg_G(z)=\delta=d-1$, so $N_G[z]=C_{i-1}\cup C_i\cup C_{i+1}$, and $z$ is exact.

Third, let $z\in C_i\cap S$ with $|i-j|\le 2$, and suppose that $z$ is not a retained vertex of positive multiplicity. If $z$ were a discarded common vertex, its used edges to $Y$ would belong to $\Omega$ and it would have no unused neighbour in $S$, so $\ell(z)=0$. If $z$ were a discarded private vertex of an unchosen cycle, its used edges to $Y$ would belong to $\Omega$, and it would have at most three unused neighbours in $S$, all on the opposite branch and of multiplicity at most one, so $\ell(z)\le 3$. If $z$ were a discarded private vertex on the positive branch of a chosen cycle, its unused neighbours in $S$ would lie on the zero branch and have multiplicity zero, so $\ell(z)=0$. All three cases contradict $\ell(z)=6$. Hence $z$ lies on the zero branch of a chosen cycle, and $a(z)=0$ whether or not $z$ is retained. The used neighbours of $z$ are not common ends of the cycle, since these are marked, so they have multiplicity zero or are joined to $z$ by edges of $\Omega$. Therefore $z$ collects its load six from at most three unused neighbours on the positive branch, and these are three vertices $a',v',b'$ of multiplicity two at the labels $\pi(z)-1$, $\pi(z)$ and $\pi(z)+1$, consecutive on the positive branch. In particular $\pi(z)=t_i$. Since $z$ does not lie on the positive path, Lemma~\ref{lem:exchange}(i) applied to $a',v',b'$ gives $\deg_G(z)\le\deg_G(v')=\delta$, so $\deg_G(z)=\delta$. There is at most one such vertex at each label. If $\delta\le4$ there is none, since its two used edges and three unused neighbours would give degree at least five.

Let $z$ be such a vertex of the zero branch with $|i-j|\le 1$. We claim that every vertex of $C_{i-1}\cup C_i\cup C_{i+1}$ other than $z$ is adjacent to $z$. Retained vertices of positive multiplicity in these sets are adjacent to $z$ by the consequence of load six. A vertex $y\in C_{i'}\setminus S$ with $|i'-i|\le 1$ is exact by the second step, as $|i'-j|\le 2$, and $z\in C_{i'-1}\cup C_{i'}\cup C_{i'+1}$, so $y$ is adjacent to $z$. By the third step, the only remaining vertices of these sets are vertices of zero branches of load six at the labels $i-1$ and $i+1$, at most one at each label. If $y$ is such a vertex at the label $i+1$, then $\pi(y)=t_{i+1}=\pi(z)+1$. Private vertices of different support cycles have labels differing by at least two, so $y$ lies on the same zero branch as $z$, and since labels strictly increase along this branch, $y$ is the successor of $z$ on it; in particular $y$ is adjacent to $z$. The label $i-1$ is symmetric. Thus $C_{i-1}\cup C_i\cup C_{i+1}\subseteq N_G[z]$. This set has $d$ vertices and $\deg_G(z)=\delta$, so $N_G[z]=C_{i-1}\cup C_i\cup C_{i+1}$, and $z$ is exact. In particular, no used edge leaves these three sets at $z$.

By the three steps, every vertex of $C_j$ is exact for every clean label $j$. Suppose that $j_0,j_0+1,\dots,j_0+20$ are clean. Then the sets $C_{j_0},\dots,C_{j_0+20}$ are nonempty and pairwise disjoint, each of them is a clique, consecutive ones are completely joined, and the conditions on the inner sets are exactly the exactness of their vertices. So they form an exact strip of length $21$, which is excluded. Hence every run of consecutive clean labels has length at most twenty. There are at most $\Nbad+1$ maximal runs of clean labels, and therefore
\[
\Lambda\le 21\Nbad+20\le 147(D_0+J+U)+735k+882s+294b+314 .
\]
If $\delta\ge4$, then $d-4\ge1$; if $\delta=3$, the hypothesis gives $k=0$. In either case, with $\kappa=1470$, we obtain
\[
2\Lambda\le\kappa\Psi+1764s+588b+628 .
\]
Every retained label carries at most two distinct support vertices, even when one of them has multiplicity zero, so $|Y|\le 2\Lambda$. Hence $\Psi\ge(|Y|-1764s-588b-628)/\kappa$, and substituting this into the bound for $T_Y/2$ above proves (ii). This completes the proof.
\end{proof}

A \emph{length-three empty window} is the set of private support vertices whose labels lie in an integer interval strictly between the labels of the two ends of one support cycle, such that no path of length at most three whose internal vertices lie outside the full support joins the two branches inside this set. The condition refers to the whole window, not only to the vertices kept in an account.

\begin{lemma}\label{lem:window}
Suppose that $G$ has minimum degree at least $\delta\ge3$ and contains no exact strip of length $21$. Let $X$ be a subset of a length-three empty window, and let $\sigma_X$ be the total number of maximal segments of consecutive path vertices of $X$ on its two branches. A \emph{full centre} is a vertex outside the full support $S$ with three neighbours in $X$. Put
\[
J_X^*=\bigcup_{z\text{ a full centre}}\bigl(N_G(z)\setminus N_G[X]\bigr),
\qquad J_X=|J_X^*|,
\qquad \eta=\frac1{2940d}.
\]
Then
\[
|X|\le\frac{3(|N_G[X]|+J_X)}d-\eta|X|
+\frac{32+(14/5)\sigma_X}{d}.
\]
Moreover, $J_X^*\subseteq N_G^2[X]\setminus N_G[X]$, and each of its vertices is adjacent to a full centre with three neighbours in $S$.
\end{lemma}

\begin{proof}
The assertion is immediate if $X=\emptyset$. Split $X$ as $X_P\cup X_Q$ according to its two branches. For each nonempty part $X_i$, choose its support cycle as the only chosen cycle, designate its branch as positive, and apply Lemma~\ref{lem:strict}(ii) with $Y=X_i$. Every retained vertex has multiplicity two, so no augmentation is needed: $k=0$ and $T_Y/2=|X_i|$. Also $s=1$ and $b=2\sigma_i$, where $\sigma_i$ is the number of path segments of $X_i$. Let $J_i^*$ be the set $J^*$ in this application, and write $h_i=|N_G[X_i]|$ and $J_i=|J_i^*|$. The lemma gives
\[
|X_i|\le\frac{3(h_i+J_i/2)}d-\eta|X_i|
+\frac{2+2392/2940+(12/5)\sigma_i}{d}.
\]
Here the centres defining $J_i^*$ are exactly the vertices outside $S$ with three neighbours in $X_i$: their three neighbours have distinct labels and account for load six.

We compare these sets for the two branches. A vertex of $N_G[X_P]\cap N_G[X_Q]$ cannot lie outside $S$, since it would give a path of length two between the branches. Nor can it lie inside the full window, since it would be incident with an edge between its two branches. A support vertex outside the window which meets both branches is incident with a used edge into the window: a private vertex has no unused neighbour on its own path, and a common vertex has no unused neighbour in $S$. Such an edge crosses one of the two label cuts at the ends of the window. Its outside end therefore belongs to one of the two sets $\partial_j$ of Lemma~\ref{lem:balls}(ii), each of order at most four. Consequently,
\[
h_P+h_Q\le |N_G[X]|+8.
\]

Every full centre has all its contacts on one branch, at the three consecutive labels $\pi(z)-1,\pi(z),\pi(z)+1$, all inside the window. Its only support neighbours are these contacts: any other support neighbour would be on the opposite branch inside the window and give a path of length two between the branches. Thus $J_P^*$ and $J_Q^*$ contain no support vertex. If a vertex of $J_P^*$ had a neighbour in $X_Q$, a contact of its defining centre, that centre, the vertex, and its neighbour in $X_Q$ would form a path of length three between the branches with internal vertices outside $S$. This is excluded, and the other direction is symmetric. Hence both $J_P^*$ and $J_Q^*$ avoid $N_G[X]$.

A full centre of $X$ lies on just one of these two sides, and the preceding observation shows that $J_X^*=J_P^*\cup J_Q^*$. In particular,
\[
\frac{J_P+J_Q}{2}\le J_X.
\]
Adding the two single-branch bounds, with an empty branch contributing zero, now proves the required estimate, since
\[
24+2\left(2+\frac{2392}{2940}\right)<32,
\qquad\frac{12}{5}\le\frac{14}{5}.
\]
The final containment follows directly from the definition of $J_X^*$.
\end{proof}

The next lemma fixes the reversals of the support cycles in a chain.

\begin{lemma}\label{lem:chain}
Suppose that some pieces form a chain in which consecutive pieces meet only in their common port. Each piece has an orientation and its complete reversal, and the prescriptions of different pieces remain compatible under every independent choice of reversals. For the $i$-th piece, let $p_i$ and $q_i$ be certified lengths of the two port-to-port routes in the two directions, let $\phi_i$ bound the directed distances in both directions between either port and every support vertex of the piece, and let $\Delta_i$ bound the directed distance between any two support vertices of the piece. Put $c_i=(p_i+q_i)/2$, $w_i=|p_i-q_i|/2$ and $u_i=\phi_i-c_i+1$. If $c_i+w_i\le\phi_i\le 2c_i-1$ for every $i$, then one fixed choice of reversals gives support diameter at most
\[
\max\Bigl\{\max_i\Delta_i,\ \sum_ic_i+u_{(1)}+u_{(2)}+u_{(3)}-3\Bigr\},
\]
where $u_{(1)}\ge u_{(2)}\ge u_{(3)}$ are the three largest values $u_i$, padded by ones if there are fewer than three pieces. Auxiliary vertices may be shared, but any shared prescribed edge must satisfy the stated compatibility.
\end{lemma}

\begin{proof}
We append the pieces one by one. For the oriented prefix, let $I$ and $O$ bound the directed distances from every support vertex of the prefix to its right port and from its right port to every support vertex of the prefix, respectively, and let $\Delta'$ bound its support diameter; for the empty prefix these three bounds are zero. We append the next piece and reverse it if necessary, so that its longer port-to-port route is added to the smaller of $I$ and $O$. Writing $p,q,\phi,\Delta$ for the parameters of the new piece in this orientation, where $p$ refers to the route from its left port to its right port and $q$ to the route back, the new prefix has the bounds
\[
I'=\max\{\phi,I+p\},\qquad O'=\max\{\phi,O+q\}\qquad\text{and}\qquad\Delta''=\max\{\Delta',\Delta,\phi+\max(I,O)\}.
\]
Let $C=\sum c_i$ over the old prefix, and let $a_1\ge a_2\ge a_3$ be its three largest values $u_i$, padded by ones. We show by induction that
\[
\frac{I+O}{2}\le C+a_1-1\qquad\text{and}\qquad\max(I,O)\le C+a_1+a_2-2 ,
\]
and that
\[
\Delta'\le\max\Bigl\{\max_i\Delta_i,\ C+a_1+a_2+a_3-3\Bigr\}.
\]
These inequalities hold for the empty prefix. Put $\theta_+=\max(I,O)$, $\theta_-=\min(I,O)$ and $C'=C+c$, and let $a_1'\ge a_2'\ge a_3'$ be the three largest values in the new list, padded by ones. Before taking maxima with $\phi$, the two new bounds are $\theta_++c-w$ and $\theta_-+c+w$. The first one is at most $C+a_1+a_2-2+c\le C'+a_1'+a_2'-2$. For the second one, $\theta_-\le(I+O)/2\le C+a_1-1$, and $w\le u-1$ since $c+w\le\phi=c+u-1$; hence $\theta_-+c+w\le C'+a_1+u-2\le C'+a_1'+a_2'-2$. The same bound holds for $\phi=c+u-1$, since $u\le a_1'$ and $a_2'\ge 1$. This proves the bound for $\max(I',O')$.

If neither of the two bounds is raised to $\phi$, their average increases by exactly $c$, so it is at most $C'+a_1-1\le C'+a_1'-1$. If both are raised, their average is $\phi=c+u-1\le C'+a_1'-1$. If exactly one is raised, their average is at most $(c+u-1+C'+a_1'+a_2'-2)/2$, and this is at most $C'+a_1'-1$ if and only if $u+a_2'-1\le C+a_1'$. The latter holds because $u\le a_1'$ and $a_2'\le C+1$: every old value satisfies $u_i\le c_i\le C$, the padding values are one, and only the new value $u$ can exceed $C$. Finally,
\[
\phi+\theta_+\le c+u-1+C+a_1+a_2-2=C'+u+a_1+a_2-3\le C'+a_1'+a_2'+a_3'-3 ,
\]
which gives the bound for $\Delta''$. This completes the induction, and the assertion follows from the bound for $\Delta'$ after the last piece. Every reversal is chosen when its piece is appended, before any ordered pair is considered. This completes the proof.
\end{proof}

A \emph{counting account} assigns a multiplicity in $\{0,1,2\}$ to every vertex of the full support, and $T$ denotes the sum of these multiplicities. Different ordered pairs may use different accounts, but the orientation is fixed before any account is chosen.

\begin{lemma}\label{lem:support}
Let $K$ be a core of a simple graph $G$, let $G_K$ be obtained from $G$ by contracting $K$ to a vertex $o$, keeping all edges between $K$ and $V(G)\setminus K$, and let $x\notin K$. Restrict $G_K$ to $o$ and the component of $G-K$ containing $x$, and choose an optimal pair of edge-disjoint $o$--$x$ paths. If the first support cycle consists of two parallel edges between $o$ and a vertex $v$, save its two edges and remove this cycle. Let $S$ be the support of the remaining pair, called the \emph{tail} when a first cycle was removed; by convention, put $S=\emptyset$ if this removal leaves no edge. When no first cycle is removed, retain exactly one edge from each parallel class, keeping the used edge whenever the class contains one. Suppose that $S\ne\emptyset$. The construction below takes place in the resulting simple quotient graph, or in $G-K$ in the tail case; neighbourhoods of its auxiliary vertices are taken there. The separation assertions in part~\textup{(v)} refer to the original graph $G$. There is one strong orientation of the subgraph induced by $S$, the selected auxiliary vertices described in the proof, and all other vertices with at least two neighbours in $S$, with the following properties.
\begin{enumerate}[label=\textup{(\roman*)}]
\item Every auxiliary vertex has directed paths to and from $S$ of length at most two.
\item For every ordered pair of vertices of $S$, there is a counting account $T$, with at most three chosen cycles and at most three accounting gaps, such that the directed distance of the pair is at most $T/2+11$.
\item Every counting account partitions $S$ as $S=Y\cup\bigcup_{X\in\mathcal{X}}X\cup S_0$, where $Y$ carries the multiplicities $w$ of this section, with the copied branches as positive branches and including its vertices of multiplicity zero, and every $X\in\mathcal{X}$ is a subset of a length-three empty window. With $x_{\mathcal{X}}=\sum_{X\in\mathcal{X}}|X|$ and $f=|S_0|$,
\[
\frac{T}{2}\le\frac{T_Y}{2}+x_{\mathcal{X}}+f\qquad\text{and}\qquad|Y|+x_{\mathcal{X}}=|S|-f .
\]
\item Let $s$ and $r$ be the numbers of chosen cycles and accounting gaps of the account, let $b$ be the number of used edges with exactly one end in $Y$, and let $\sigma$ be the total number of maximal path segments of the sets $X\in\mathcal{X}$. Then $s\le3$, $r\le3$, $|\mathcal{X}|\le2r\le6$, and
\[
f\le688r+280\le2344,\qquad b\le8r+4f\le9400\qquad\text{and}\qquad\sigma\le4r+f\le2356 .
\]
\item The sets $N_G^2[Y]$ and $N_G^2[X]$, $X\in\mathcal{X}$, are pairwise disjoint and disjoint from $N_G[K]$.
\end{enumerate}
\end{lemma}

\begin{proof}
By the convention in the statement, a nonempty support has an edge.

\emph{Rungs.} Consider a support cycle with ends $A$ and $B$, and write its two branches as $P_0=A,P_1,\dots,P_a=B$ and $Q_0=A,Q_1,\dots,Q_b=B$, in increasing order of labels. A \emph{rung} is a path of length at most three between a private vertex of one branch and a private vertex of the other branch whose internal vertices lie outside the full support; all its edges are unused. The \emph{interval} of a rung is the smallest integer interval containing the labels of all its vertices, including the internal ones. Since all edges of a rung are unused, this interval has width at most three, and since every internal vertex of a rung is adjacent to one of its ends, the interval lies in the label interval of the cycle. We choose a family of rungs, over all support cycles at once, whose intervals are pairwise disjoint and which is maximal under inclusion among such families. In each cycle, the selected rungs are ordered along the cycle, and both their ends on $P$ and their ends on $Q$ strictly increase along this order, since labels strictly increase along each branch. We add two formal rungs of length zero, one joining $P_0$ and $Q_0$ and one joining $P_a$ and $Q_b$. We direct $P$ forwards and $Q$ backwards, and we direct the selected rungs alternately from $P$ to $Q$, called \emph{positive} rungs, and from $Q$ to $P$, called \emph{negative} rungs, along the order; each formal rung serves as a rung of both types (see Figure~\ref{fig:rungs}). Thus between two consecutive rungs of one type there is at most one rung of the other type. These prescriptions are compatible with each other and with independent reversals of whole cycles, since selected rungs have disjoint intervals and hence disjoint internal vertices, their ends are distinct private vertices, and no rung edge is used by the pair.

\begin{figure}[ht]
\centering
\begin{tikzpicture}[>=Stealth, vx/.style={circle,fill,inner sep=1.3pt}, ax/.style={circle,draw,fill=white,inner sep=1.4pt}]
\node[vx,label=left:$A$] (A) at (0,0) {};
\node[vx,label=right:$B$] (B) at (11,0) {};
\foreach \i in {1,...,7} {\node[vx] (P\i) at (1.375*\i,1.3) {};}
\foreach \i in {1,...,6} {\node[vx] (Q\i) at (1.5714*\i,-1.3) {};}
\draw[->] (A) -- (P1);
\foreach \i/\j in {1/2,2/3,3/4,4/5,5/6,6/7} {\draw[->] (P\i) -- (P\j);}
\draw[->] (P7) -- (B);
\draw[->] (B) -- (Q6);
\foreach \i/\j in {6/5,5/4,4/3,3/2,2/1} {\draw[->] (Q\i) -- (Q\j);}
\draw[->] (Q1) -- (A);
\draw[->,very thick] (P2) -- node[right] {$+$} (Q2);
\node[ax] (Y1) at (6.0,0) {};
\draw[->,very thick] (Q4) -- (Y1);
\draw[->,very thick] (Y1) -- node[right] {$-$} (P4);
\draw[->,very thick] (P6) -- node[left] {$+$} (Q5);
\node at (5.5,1.85) {$P$};
\node at (5.5,-1.85) {$Q$};
\end{tikzpicture}
\caption{A support cycle with ends $A$ and $B$. The branch $P$ is directed from $A$ to $B$, the branch $Q$ from $B$ to $A$, and the selected rungs (thick) are directed alternately from $P$ to $Q$ and from $Q$ to $P$. The white vertex is an internal vertex of a rung of length two; the formal rungs at $A$ and $B$ have length zero.}
\label{fig:rungs}
\end{figure}
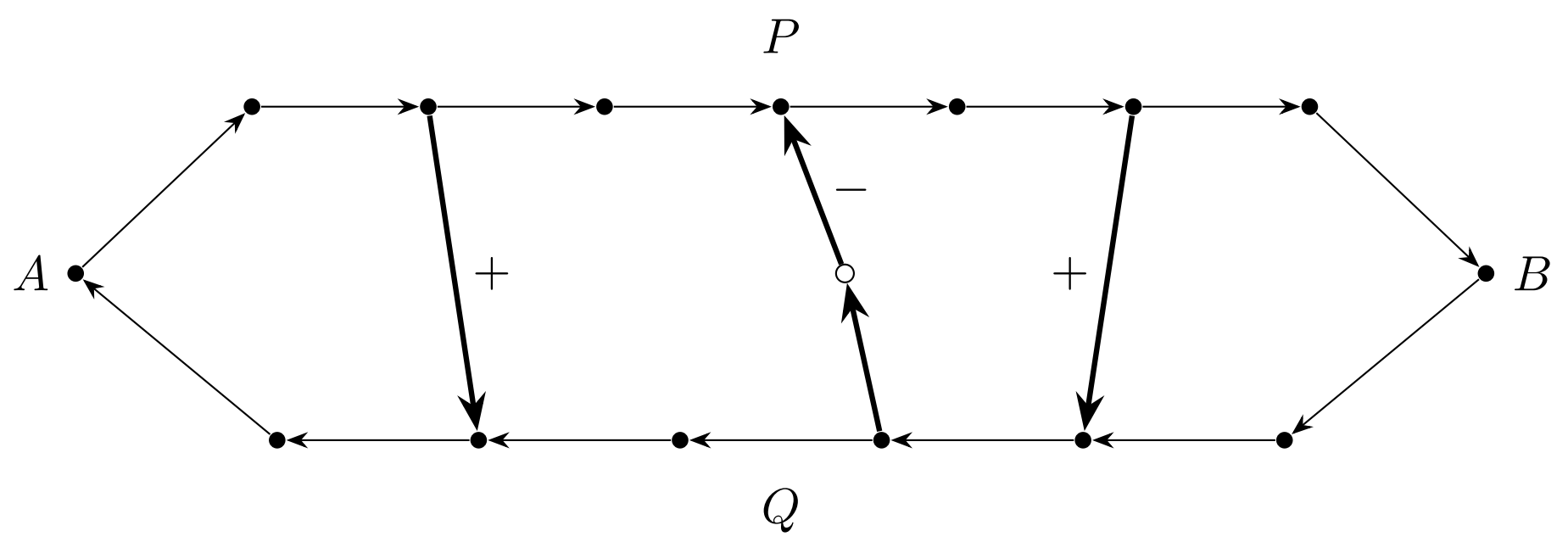

For a rung $\varrho$ we write $\varrho^P$ and $\varrho^Q$ for the indices of its ends on $P$ and $Q$, and $|\varrho|$ for its length. Let $g_P$ be the largest difference $\varrho'^P-\varrho^P$ over all pairs of consecutive rungs $\varrho,\varrho'$ of the same type, including the formal rungs, and define $g_Q$ in the same way. With $l=a+b$ and $c=l/2$ put
\[
\phi=\min\{\,l-1,\ \max(a+g_Q,\ b+g_P)+3\,\}.
\]
We claim that $\phi$ bounds the directed distances in both directions between each end of the cycle and each vertex of the cycle. The forward branch gives distance at most $a$ from $A$ to $P_i$ and from $P_i$ to $B$, and the backward branch gives distance at most $b$ from $B$ to $Q_i$ and from $Q_i$ to $A$. To go from $P_i$ to $A$, follow $P$ to the next positive rung $\varrho$, so that $\varrho^P-i\le g_P$, cross it, and follow $Q$ to $A$; this costs at most $g_P+3+b$. To go from $B$ to $P_i$, follow $Q$ to the last negative rung $\varrho$ with $\varrho^P\le i$, cross it, and follow $P$ to $P_i$; again this costs at most $b+3+g_P$. Symmetrically, the routes from $A$ to $Q_i$ and from $Q_i$ to $B$ cost at most $a+3+g_Q$. Finally, the directed cycle formed by $P$ and $Q$ gives the alternative bound $l-1$. Hence $\max(a,b)\le\phi\le l-1$, as required in Lemma~\ref{lem:chain} with port-to-port lengths $a$ and $b$.

\emph{Requests inside one cycle.} A request in the direction of a branch uses that branch. For a request from $P_i$ to $P_j$ with $i>j$, let $\varrho$ be the first positive rung with $\varrho^P\ge i$ and $\varrho'$ the last negative rung with $\varrho'^P\le j$. Since $\varrho'$ precedes $\varrho$, the route along $P$ from $P_i$ to $P_{\varrho^P}$, across $\varrho$, along $Q$ to $Q_{\varrho'^Q}$, across $\varrho'$ and along $P$ to $P_j$ has length
\[
(\varrho^P-i)+|\varrho|+(\varrho^Q-\varrho'^Q)+|\varrho'|+(j-\varrho'^P)\le b+(\varrho^P-i)+(j-\varrho'^P)+6 .
\]
The two segments of $P$ used here, called \emph{tails}, are vertex-disjoint since $j<i$, and each lies in a gap between consecutive rungs of one type. A request from $Q_i$ to $Q_j$ with $i<j$ is symmetric and costs at most $a+(i-\varrho'^Q)+(\varrho^Q-j)+6$, with two disjoint tails on $Q$. For a request from $P_i$ to $Q_j$, let $\varrho$ be the first positive rung with $\varrho^P\ge i$ and $\varrho^Q\ge j$, and use the route along $P$ to $P_{\varrho^P}$, across $\varrho$ and along $Q$ to $Q_j$. Unless the request starts at $A$, which is trivial, the preceding positive rung $\varrho_-$ fails one of the two comparisons. If $\varrho_-^P<i$, the segment of $P$ from $P_i$ to $P_{\varrho^P}$ is a tail in the gap between $\varrho_-$ and $\varrho$, and the rest of the route costs at most $b+3$. If $\varrho_-^Q<j$, the segment of $Q$ from $Q_{\varrho^Q}$ to $Q_j$ is a tail in that gap, and the rest costs at most $a+3$. A request from $Q_j$ to $P_i$ uses the last negative rung satisfying both reverse comparisons, symmetrically. Thus every request inside one cycle costs at most the length of one whole branch, plus the total length of one or two disjoint tails on the opposite branch, each inside a gap between consecutive rungs of one type, plus at most six further edges. We apply Lemma~\ref{lem:chain} once, with the pieces being the support cycles, to fix the reversal of every cycle; this happens before any account is chosen.

\emph{Gaps and windows.} An \emph{accounting gap} runs from the lower end of the interval of a selected rung to the upper end of the interval of the next selected rung of the same type; the intervals of the formal rungs are single labels. For a union of $r$ gaps, we mark the two end intervals of each gap and the at most one interval of the opposite type in between. This gives at most $3r$ marked intervals, which leave at most $2r$ complementary open integer intervals inside the gaps; taking the union avoids double counting when gaps overlap. Each marked interval contains at most eight support vertices, namely at most two at each of four labels. In each complementary integer interval $[l',u']$ we also discard the two end labels and keep $[l'+1,u'-1]$ if it is nonempty; this discards at most $8r$ further vertices. The kept set of private vertices is a length-three empty window. Indeed, a rung with both ends in $[l'+1,u'-1]$ has its whole interval in $[l',u']$, since every internal vertex of a rung is adjacent to one of its ends. Such a rung would be disjoint from every selected interval: inside the gaps the selected intervals are exactly the marked ones, because between two consecutive rungs of one type there is at most one rung of the other type, and selected intervals of other cycles lie in other label ranges. This contradicts the maximality of the selected family.

At both ends of every kept window we apply the deletion of Lemma~\ref{lem:balls}(ii) at the corresponding cut; there are at most $4r$ cuts, and each deletion discards at most $164$ support vertices. We also discard every support vertex within distance three of $o$; by Lemma~\ref{lem:balls}(iii) there are at most $280$ of them. Let $S_0$ be the set of all discarded vertices. Then
\[
f=|S_0|\le24r+8r+4r\cdot164+280=688r+280 .
\]
The vertices of $S\setminus S_0$ outside the kept windows form the fragment $Y$, which may consist of several pieces, and the kept windows minus $S_0$ form the sets $X\in\mathcal{X}$. A vertex of $Y$ and a vertex of a set $X\in\mathcal{X}$, as well as two vertices of different sets of $\mathcal{X}$, lie on opposite sides of a cut at which the deletion of Lemma~\ref{lem:balls}(ii) was applied; recall that all support vertices with labels in the interval of a kept window belong to that window or to $S_0$. By Lemma~\ref{lem:balls}(ii), the remark following it, and Lemma~\ref{lem:balls}(iii), this proves (v). Initially the windows have at most $4r$ segments, two for each of at most $2r$ windows, and each discarded vertex adds at most one more. Each path crosses a cut at most once, and a discarded vertex has at most four used edges. This proves the three bounds in (iv); the bounds $s\le 3$ and $r\le 3$ follow from the choice of the accounts below.

\emph{Accounts.} Initially give multiplicity two to every common vertex and one to every private vertex. \emph{Copying} a branch of a cycle gives multiplicity two to its private vertices and zero to the private vertices of the opposite branch. In a selected gap, both private multiplicities are restored to two. Then every vertex of a kept window has multiplicity two, and the fragment $Y$ carries exactly the multiplicities $w$ with the copied branches as positive branches; indeed, every vertex of a restored gap lies in a kept window or in $S_0$. Every multiplicity is at most two, which proves (iii).

Let $\mu$ be the total length of the pair, after the removal of a parallel first cycle and the simplification; it has $\mu+2$ vertex occurrences. For the top-three term of Lemma~\ref{lem:chain}, choose its $q\le 3$ contributing cycles. In each of them, choose a branch of length $L_i$ and an opposite gap of length $g_i$ attaining the maximum in the definition of $\phi_i$, and let $L_i'$ be the length of the other branch. Copy the branch of length $L_i$ and restore its gap. Since the opposite subpath inside the gap has at least $g_i-1$ private vertices,
\[
T\ge\mu+2+\sum_{i=1}^q(L_i-L_i')+2\sum_{i=1}^q(g_i-1).
\]
Since $\phi_i\le L_i+g_i+3$ and $c_i=(L_i+L_i')/2$, the top-three term of Lemma~\ref{lem:chain}, with the padding values included, is at most
\[
\frac{\mu}{2}+\sum_{i=1}^q(\phi_i-c_i)\le\frac{T}{2}+4q-1\le\frac{T}{2}+11 .
\]
For a request inside one cycle, let $L$ be the length of the whole branch which it uses, and let $t$ be the total length of its one or two tails on the opposite branch. Copy that branch and restore the union of the gaps containing the tails. The tails are disjoint and contain at least $t-2$ distinct private vertices. The copied branch contributes at least $2L+2$, namely $2(L-1)$ from its private vertices and four from its ends, so $T\ge 2L+2+2(t-2)$ and $T/2\ge L+t-1$. Hence the route has length at most $L+t+6\le T/2+7$. A request along one branch needs no gap. These accounts have $s=1$ and $r\le 2$, while the top-three accounts have $s=q\le 3$ and $r=q\le 3$. Together with Lemma~\ref{lem:chain}, this proves (ii), and also completes the proof of (iv).

\emph{Auxiliary vertices.} Finally, we add the internal vertices of all selected rungs and every other vertex with at least two neighbours in $S$. A selected rung gives each of its internal vertices directed paths to and from $S$ of length at most two. Every other such vertex lies on no selected rung, so none of its edges has been prescribed; we orient one of its edges to $S$ towards it and another one, to a different vertex of $S$, away from it. All remaining edges are oriented arbitrarily. This gives one strong orientation in which every auxiliary vertex has access to and from $S$ within two steps, which proves (i). This completes the proof.
\end{proof}

For the rest of the proof, use the uniform constants
\[
\kappa=1470,\qquad \eta=\frac1{2940d},\qquad
(a,b_0,c)=\left(\frac{13}{5},\frac65,\frac{157}{735}\right).
\]
The positive terms of Lemma~\ref{lem:strict}(ii) are
\[
\frac{1764s+588b+628}{2940d}+\frac{2s+b}{d}
=\frac{as+b_0b+c}{d}.
\]

\begin{proposition}\label{prop:extension}
Let $G$ be a bridgeless graph of minimum degree at least $\delta\ge4$ containing no exact strip of length $21$. Let $K$ be a core with prescribed orientation of diameter $\Delta_K$, and let $x\notin K$. Choose the optimal pair after contracting $K$, and define its support or tail $S$ as in Lemma~\ref{lem:support}, including the empty-tail convention. Then there are a set $B\supseteq K\cup\{x\}$ and a strong orientation $D_B$ of $G[B]$ extending the prescribed orientation such that, with $m=|N_G[B]\setminus N_G[K]|$,
\[
\diam(D_B)\le\Delta_K+\frac{3m}{d}-\eta|S|+6000.
\]
\end{proposition}

\begin{proof}
If the first support cycle consists of two parallel edges and the tail is empty, then $x$ has two neighbours in $K$. We put $B=K\cup\{x\}$ and orient one of the two edges towards $x$ and the other one away from $x$. Then $\diam(D_B)\le\Delta_K+2$, so the assertion holds since $S=\emptyset$.

Assume now that $S\ne\emptyset$. We restrict the contracted graph to the root $o$ and the component of $G-K$ containing $x$. This multigraph is connected and bridgeless, since an edge whose removal disconnects it would also disconnect $G$. We keep every edge between $K$ and this component, with its original end in $K$, choose the optimal pair with a fixed score for the root, and take the orientation and the accounts of Lemma~\ref{lem:support}.

Every retained vertex of an account has distance at least four from $K$. Hence its second neighbourhood avoids $N_G[K]$, every vertex whose degree is used in Lemmas~\ref{lem:strict} and~\ref{lem:window} has its original degree and its original incident edges, and all exchanges in these lemmas avoid the root. The three priorities of the pair survive the simplification of parallel edges, since an improvement there would give an improvement of the original pair; if a parallel first cycle was removed, the same holds after attaching its two saved edges. Hence Lemmas~\ref{lem:strict} and~\ref{lem:window} apply, and every exact strip which they could produce would be an exact strip of $G$.

Let $B$ consist of $K$, the support $S$ and the auxiliary vertices of Lemma~\ref{lem:support}. The vertices of $J^*$ in Lemma~\ref{lem:strict} are neighbours of vertices outside the support with three neighbours of positive multiplicity in $S$, and the vertices of $J_X^*$ in Lemma~\ref{lem:window} are neighbours of full centres. All these centres have at least two neighbours in $S$, so they are auxiliary vertices and belong to $B$. Consequently $J^*$ and all sets $J_X^*$ are contained in $N_G[B]$. Since $J^*\subseteq N_G^2[Y]$ and $J_X^*\subseteq N_G^2[X]$, Lemma~\ref{lem:support}(v) shows that the sets $N_G[Y]\cup J^*$ and $N_G[X]\cup J_X^*$, $X\in\mathcal{X}$, are pairwise disjoint subsets of $N_G[B]\setminus N_G[K]$. Writing $h_Y=|N_G[Y]|$, $J_Y=|J^*|$ and $h_X=|N_G[X]|$, we obtain
\[
h_Y+J_Y+\sum_{X\in\mathcal{X}}(h_X+J_X)\le m .
\]
This counts distinct vertices of $G$; the vertices of $J^*$ and of the sets $J_X^*$ need not belong to $B$.

Now fix any account, and put $w_{\mathcal{X}}=|\mathcal{X}|$. By Lemma~\ref{lem:support}(iii), Lemma~\ref{lem:strict}(ii) for $Y$ (with its $J/2$ term relaxed to $J$), Lemma~\ref{lem:window} for every $X\in\mathcal{X}$, and the last inequality,
\[
\frac{T}{2}\le\frac{T_Y}{2}+x_{\mathcal{X}}+f
\le\frac{3m}{d}-\eta|S|+(1+\eta)f
+\frac{as+b_0b+c+32w_{\mathcal{X}}+(14/5)\sigma}{d}.
\]
The negative term uses $|Y|+x_{\mathcal{X}}=|S|-f$, including the retained vertices of multiplicity zero. With $s\le3$, $f\le2344$, $b\le9400$, $\sigma\le2356$ and $w_{\mathcal{X}}\le6$, the positive terms, together with eleven steps for support routing, four steps for auxiliary ends and one possible saved edge, are at most
\[
(1+\eta)2344
+\frac{(13/5)3+(6/5)9400+157/735+192+(14/5)2356}{d}+16
<5976<6000.
\]
The left-hand side is decreasing for $d\ge5$ and is less than $5975.523$ at $d=5$.

Suppose first that there is no parallel first cycle. For each ordered pair of quotient vertices, Lemma~\ref{lem:support}(i),(ii) gives a route of length at most $T/2+15$ for a suitable account. The preceding uniform estimate bounds the quotient diameter, and Lemma~\ref{lem:lift} lifts this orientation while preserving $K$, at an additional cost of at most $\Delta_K$. This gives the required bound for pairs not both in $K$. Unused edges are oriented arbitrarily.

Suppose next that the first support cycle consists of two parallel edges $e_1,e_2$ between $v$ and $K$. Since $G$ is simple, $e_1$ and $e_2$ have distinct ends in $K$. We apply the construction of Lemma~\ref{lem:support} to the tail in $G-K$, keeping the labels of the original pair and the deletions near the root. Every retained vertex still has distance at least four from $K$, so a path through $K$ between two retained vertices has length at least eight, and the separation in Lemma~\ref{lem:support}(v) also holds in $G$. The vertex $v$ is the only vertex of the tail support which can be joined to $K$ by an unused edge: an unused edge from $o$ to a later support vertex $y$ would give $\pi(y)\le 1$, whereas $\pi(y)\ge\pi(v)+1\ge 2$. We orient $e_1$ from $K$ to $v$ and $e_2$ from $v$ to $K$. A request between two vertices of $B\setminus K$ uses the orientation of the tail. A request between a vertex of $K$ and a vertex of $B\setminus K$ uses one route inside $K$, one of the edges $e_1,e_2$, and a route of the tail from or to $v$. Hence the bound increases by at most one, and it remains below $\Delta_K+3m/d-\eta|S|+6000$.

Pairs of vertices of $K$ keep their old routes. To check that their bound is also valid, take any account in the preceding estimate. Since $T\ge0$, it implies $3m/d-\eta|S|+6000\ge0$. Thus the old core diameter is not larger than the asserted bound. This completes the proof.
\end{proof}

\section{Minimum degree three and the initial core}\label{sec:small}

The window estimate in Lemma~\ref{lem:window} applies also when $\delta=3$. The remaining fragment needs a separate count, because the coefficient $d-4$ in Lemma~\ref{lem:strict} vanishes. All degrees and neighbourhoods below refer to the fixed simple graph before contraction, and the pairs have the same three priorities as before.

\begin{lemma}\label{lem:fragment3}
Let $G$ have minimum degree at least three, and let $Y$, $w$, $T_Y$, $h$, $s$ and $b$ be as in Section~\ref{sec:counting}. Let $\nu$ be the number of common vertices of the full pair other than its two ends, let $U_1$ be the number of retained private vertices in unchosen cycles, and let $z_1$ be the number of retained private vertices of multiplicity zero in chosen cycles. Suppose that at most $f$ common vertices other than the two ends of the pair do not belong to $Y$. Then
\[
\frac{T_Y}{2}\le\frac{3h}{4}-\frac{\nu}{8}-\frac{U_1+z_1}{32}+\frac{f}{4}+\frac{s}{2}+\frac{7b}{32}+\frac14 .
\]
\end{lemma}

\begin{proof}
Let $\varepsilon_u$, $\varepsilon_+$ and $\varepsilon_c$ be the sums of $\deg_G(v)-3$ over the retained private vertices of unchosen cycles, the retained private vertices of positive branches and the retained common vertices, respectively. Then
\[
W=\sum_{v\in Y}w(v)\bigl(\deg_G(v)+1\bigr)=4T_Y+\varepsilon_u+2\varepsilon_++2\varepsilon_c .
\]
Every common vertex other than the two ends of the pair has four distinct used edges and hence degree at least four, so $\varepsilon_c\ge\nu-f$.

For every centre $z\in N_G[Y]$ let $L(z)=\sum_{v\in N_G[z]\cap Y}w(v)$ be its closed load with respect to the multiplicities $w$; these loads sum to $W$. A centre outside the support has load at most six. A common centre has load at most ten, and there are at most $\nu+2$ common centres, which contribute at most $4(\nu+2)$ above six. For a private centre $v$ of an unchosen cycle, let $\alpha$ be the number of its used neighbours which are common. If $v$ is retained, its load is at most $\deg_G(v)+1+\alpha$; if it is discarded, its load is at most $5+\alpha$, since it has at most three unused neighbours in the support and two used neighbours. There are at most $4(\nu+2)$ used incidences between common and private vertices. Summing, these centres contribute at most $-2U_1+\varepsilon_u+4(\nu+2)$ above six.

A private centre on a positive branch has load at most six, even if it is discarded. For a centre $z$ on a zero branch, let $k'$ be the number of its retained private neighbours on the positive branch, and let $\alpha$ again be the number of its used common neighbours. Then $k'\le 3$ and $L(z)\le 2k'+2\alpha$, so a retained such centre has load at most $4+2\alpha$ unless $k'=3$. In the exceptional case $k'=3$, the three positive neighbours of $z$ are consecutive on the other branch, and the two used edges of $z$ give $\deg_G(z)\ge 5$. The exchange of Lemma~\ref{lem:exchange}(i) then forces the middle positive neighbour to have degree at least five as well, and we charge $z$ to this neighbour through their unused edge. A positive vertex $w$ receives at most $\deg_G(w)-2\le 3(\deg_G(w)-3)/2$ charges, so twice the number of exceptional centres is at most $3\varepsilon_+$. A discarded centre on a zero branch has load at most $6+2\alpha$. In every chosen cycle there are at most two used incidences between private vertices of the zero branch and common vertices. Hence the centres on zero branches contribute at most $-2z_1+3\varepsilon_++4s$ above six.

Combining all cases gives $W\le 6h+8(\nu+2)+4s-2(U_1+z_1)+\varepsilon_u+3\varepsilon_+$, and hence
\[
\frac{T_Y}{2}\le\frac{3h}{4}+\nu+2+\frac{s}{2}-\frac{U_1+z_1}{4}+\frac{\varepsilon_+}{8}-\frac{\varepsilon_c}{4}.
\]
On the other hand, Lemma~\ref{lem:strict}(i), which is valid at degree three, gives
\[
\frac{T_Y}{2}\le\frac{3h}{4}+\frac{s}{2}+\frac{b}{4}-\frac{\varepsilon_u}{8}-\frac{\varepsilon_+}{4}-\frac{\varepsilon_c}{4}.
\]
We take one eighth of the first bound and seven eighths of the second. The coefficients of $\varepsilon_u$ and $\varepsilon_+$ become $-7/64$ and $-13/64$, and $\nu/8-\varepsilon_c/4\le-\nu/8+f/4$. Dropping the two remaining nonpositive terms proves the lemma. This completes the proof.
\end{proof}

\begin{proposition}\label{prop:extension3}
Let $\delta=3$, so that $d=4$, and let $G$ be a bridgeless graph of minimum degree at least three containing no exact strip of length $21$. Let $K$ be a core whose prescribed orientation has diameter $\Delta_K$, and suppose that $R=\max_v\dist_G(v,K)>0$. Then there are a set $B\supseteq K$ containing a vertex farthest from $K$ and a strong orientation $D_B$ of $G[B]$ extending the prescribed orientation such that, with $m=|N_G[B]\setminus N_G[K]|$,
\[
\diam(D_B)\le\Delta_K+\frac{3m}{4}-\frac{R}{11760}+7000.
\]
\end{proposition}

\begin{proof}
Put $\eta_3=1/11760$. Choose a farthest vertex $x$ and an optimal pair after contracting $K$. Use the support orientation of Lemma~\ref{lem:support}, and let $B$ consist of $K$, its support or tail, and its auxiliary vertices, as in the proof of Proposition~\ref{prop:extension}. If the pair consists only of a parallel first cycle, then $R=1$ and two oppositely oriented edges attach $x$ at a cost of at most two, so the assertion follows. We henceforth assume that the support or tail is nonempty.

First suppose there is no parallel first cycle. For an account, let $Y$ be its fragment, let $x_{\mathcal X}$ be the total number of retained window vertices, let $w=|\mathcal X|$, and let $\sigma$ be the total number of their path segments. Write $h=|N_G[Y]|$ and $h_X=|N_G[X]|$. The sets $J_X^*$ supplied by Lemma~\ref{lem:window} lie in $N_G[B]$, since their centres have three neighbours in the support and hence belong to $B$. By the separation in Lemma~\ref{lem:support}(v),
\[
h+\sum_{X\in\mathcal X}(h_X+J_X)\le m,
\qquad
\frac T2\le\frac{T_Y}{2}+x_{\mathcal X}+f.
\]
Every common vertex other than the two ends which is not in $Y$ belongs to $S_0$, because windows contain only private vertices.

Choose the zero branch in each chosen cycle and either branch in each unchosen cycle. Their concatenation is a path from the contracted root to $x$. Its internal common vertices are the $\nu$ vertices counted in Lemma~\ref{lem:fragment3}; every private vertex on it is counted by $U_1$, by $z_1$, by a retained window, or by $S_0$. Thus
\[
R\le\nu+1+U_1+z_1+x_{\mathcal X}+f.
\]
Apply Lemma~\ref{lem:fragment3} to $Y$ and Lemma~\ref{lem:window}, with $d=4$, to every window. Since $1/8,1/32\ge\eta_3$, the preceding two inequalities give
\[
\frac T2\le\frac{3m}{4}-\eta_3R
+\frac{5f}{4}+\eta_3(f+1)+\frac s2+\frac{7b}{32}+\frac14
+\frac{32w+(14/5)\sigma}{4}.
\]
Here $s\le3$, $w\le2r\le6$, $f\le2344$, $b\le9400$ and $\sigma\le2356$ by Lemma~\ref{lem:support}. These positive terms, plus eleven steps for support routing, four for auxiliary ends, and the extra $1+\eta_3$ used below for a parallel first cycle, are less than $6702$, and hence less than $7000$.

If the first support cycle is parallel, remove it and apply the same calculation to the tail. Its three priorities are inherited by reattaching the two saved edges. A tail path to $x$ has length at least $R-1$, so the last estimate loses at most $\eta_3$; a route between $K$ and the tail uses at most one saved edge. This is the additional $1+\eta_3$ already allowed above. The deletion near the original contracted root ensures that the degrees, exchanges and disjoint neighbourhood charges are interpreted in $G$, exactly as in Proposition~\ref{prop:extension}.

Lift the quotient routes through $K$, paying at most one $\Delta_K$, by Lemma~\ref{lem:lift}; in the parallel case use the saved incoming or outgoing edge as appropriate. Pairs inside $K$ keep their original routes. To check their bound as well, note that a used root--$x$ path contains a vertex at each distance $2,\ldots,R$ from $K$, so $m\ge R-1$. Hence $3m/4-\eta_3R+7000>0$. Orienting all remaining edges arbitrarily completes the extension.
\end{proof}

For a connected set $K$ and a vertex $x\notin K$, let $\mu_K(x)$ be the minimum total length of two edge-disjoint paths from the contracted root to $x$ in $G/K$, keeping parallel edges. This quantity is finite in a bridgeless graph. A single stopping condition now gives the initial core for every $\delta\ge4$.

\begin{proposition}\label{prop:initial}
Let $G$ be a bridgeless graph of minimum degree at least $\delta\ge4$ containing no exact strip of length $21$, and put $d=\delta+1$. Then $G$ has a core $K_0$ such that
\[
\diam(D_{K_0})\le\frac{3|N_G[K_0]|}{d},
\qquad
\mu_{K_0}(x)\le36\,000\,000d\quad(x\notin K_0).
\]
\end{proposition}

\begin{proof}
Start with a one-vertex core. While some outside vertex $x$ satisfies $\mu_K(x)>36\,000\,000d$, choose its optimal pair and apply Proposition~\ref{prop:extension}. The pair has $\mu_K(x)+2$ vertex occurrences, each vertex occurs at most twice, and deleting a parallel first cycle loses at most the contracted root from its support. In particular,
\[
\mu_K(x)\le2|S|+2,
\qquad
\eta|S|>\frac{18\,000\,000d-1}{2940d}>6000.
\]
Thus each extension increases the diameter by at most $3|N_G[B]\setminus N_G[K]|/d$. It adds $x$, so the process terminates. The charged sets are differences of nested closed neighbourhoods in the same graph; their sizes sum to at most $|N_G[K_0]|$. This proves the diameter bound, and the stopping rule proves the two-path bound.
\end{proof}

\begin{corollary}\label{prop:initialsmall}
Let $\delta=3$, so that $d=4$, and let $G$ be a bridgeless graph of minimum degree at least three containing no exact strip of length $21$. Then $G$ has a core $K_0$ satisfying
\[
\diam(D_{K_0})\le\frac{3|N_G[K_0]|}{4},
\qquad
\max_v\dist_G(v,K_0)\le82\,320\,000.
\]
\end{corollary}

\begin{proof}
Start with a one-vertex core and apply Proposition~\ref{prop:extension3} towards a farthest vertex whenever its distance exceeds $82\,320\,000=7000\cdot11760$. Each extension then costs at most three quarters of the number of newly charged vertices. The same telescoping argument as above proves the diameter bound, and termination gives the radius bound.
\end{proof}

\section{A density bound preserving a prescribed core}\label{sec:density}

We need to finish an orientation when few vertices remain outside the closed neighbourhood of the prescribed core. The argument of this section does not use the counting lemmas of the previous sections. We first prove a bound for an ordinary graph: a sparse fixed subgraph guarantees strong connectivity, while short ears leave independent choices which are used to shorten routes.

\begin{lemma}\label{lem:sparse}
Let $G$ be a bridgeless graph of order $N$ and minimum degree at least $d-1$, where $d\ge 3$. Then $G$ has a bridgeless subgraph $F$, possibly a single vertex, such that every vertex of $G$ is within distance two of $V(F)$ and $|E(F)|\le 18N/d$.
\end{lemma}

\begin{proof}
Start with a tree $T$ consisting of one vertex. While some vertex of $G$ has distance at least three from $T$, choose a vertex at distance exactly three from $T$, which exists by connectivity, and add a shortest path of length three joining it to $T$. The chosen vertices, together with the initial vertex, are pairwise at distance at least three, so their closed neighbourhoods are disjoint and each of them has at least $d$ vertices. Hence there are at most $N/d$ additions, and $|E(T)|\le 3N/d$.

Extend $T$ to a spanning tree $T^+$ of $G$ by a breadth-first search starting simultaneously at all vertices of $T$. Every vertex outside $T$ reaches its first ancestor in $T$ in at most two tree edges. For every $e\in E(T)$, the cut of $T^+-e$ contains an edge $g\ne e$ of $G$, since $e$ is not a bridge of $G$, and we add the fundamental cycle of $g$ with respect to $T^+$. Outside $T$, this cycle uses at most five edges, namely $g$ and the two paths of length at most two from the ends of $g$ to $T$.

Let $F$ be the union of $T$ and all these cycles. It is connected, and each of its edges lies on one of the added cycles, since the fundamental cycle added for $e\in E(T)$ passes through $e$. Thus $F$ is bridgeless and $|E(F)|\le 6|E(T)|\le 18N/d$. Its vertex set contains $V(T)$, so every vertex of $G$ is within distance two of $V(F)$. If $T$ has no edge, we take $F=T$. This completes the proof.
\end{proof}

The next lemma provides independent choices which never destroy strong connectivity. Here an \emph{ear} is a path or a closed path whose internal vertices are distinct and lie outside the old subgraph and whose ends lie in the old subgraph; if its ends coincide, it is a closed ear. Physical edges are always distinct.

\begin{lemma}\label{lem:ears}
Let $G$ be a bridgeless graph of minimum degree at least two, and let $F$ be a connected subgraph of $G$ with a strong orientation such that every vertex of $G$ is within distance two of $B=V(F)$. Then the vertices outside $B$ can be added in at most five rounds of ears, each of length at most five, such that ears of the same round have disjoint interiors and all ears are edge-disjoint. For every choice of a direction of each entire ear, the union of the ears with $F$ is strong. Moreover, every vertex has a directed path to $B$ and a directed path from $B$, each of length at most $25$, which can be found by inspecting the directions of at most five ears.
\end{lemma}

\begin{proof}
We contract $B$ to a root $o$, only in order to choose the ears, keeping all parallel edges together with their actual ends in $B$. Every neighbour of $o$ is called a \emph{top}. Every other vertex outside $B$ has a neighbour which is a top, and we fix one such neighbour as its \emph{parent}. A top together with the vertices whose parent it is forms a \emph{branch}; the branches are disjoint and cover $V(G)\setminus B$. Two branches are \emph{adjacent} if an edge of $G$ joins them, and for every pair of adjacent branches we fix one such edge as its \emph{representative edge}.

\emph{First round.} A branch which is adjacent to no other branch has at least two edges to $o$, since otherwise its only edge to the rest of the graph would be a bridge. All these edges are incident with its top, because the other vertices of the branch are not adjacent to $B$. We use two of them as an ear of length two which adds the top. For the remaining branches, we choose a maximal matching in the graph whose vertices are these branches and whose edges are the adjacent pairs. For each matched pair of branches with tops $t$ and $t'$, we take an edge from $t$ to $B$, the representative edge between the two branches, an edge from $t'$ to $B$, and, if needed, the edges from the ends of the representative edge to their parents. This is an ear of length at most five, and different matched pairs use disjoint branches.

\emph{Second and third rounds.} Every top $t$ which has not been added yet belongs to an unmatched branch which is adjacent to some other branch. By the maximality of the matching, the branch of $t$ is adjacent to a matched branch, say with top $t'$. We form an ear which starts with an edge from $B$ to $t$, possibly continues to a vertex of the branch of $t$, and then uses the representative edge to the branch of $t'$. If the other end of the representative edge has already been added, the ear stops there. Otherwise this end is a vertex $z$ of the branch of $t'$, and the ear continues through $z$ to $t'$, which was added in the first round. These ears have length at most four, and the only possible conflict between them is that several of them use the same new vertex $z$. We group the ears according to this vertex $z$. In the second round we add one ear from each group, together with all ears whose representative edge ends at a vertex already added; their interiors are disjoint. In the third round we add the remaining ears, which now stop at the vertex $z$ added in the second round. Their branches are distinct, so their interiors are again disjoint. After the third round all tops have been added.

\emph{Fourth and fifth rounds.} Every remaining vertex has a neighbour already added, namely its parent. In the subgraph induced by the remaining vertices, we take a maximal matching. Each matched edge, together with one edge from each of its ends to a vertex already added, forms an ear of length three, and each isolated vertex of this subgraph has at least two neighbours already added and forms an ear of length two with two of its edges. We add all these ears in the fourth round. Every vertex that is still left is unmatched and has a neighbour in the subgraph, and by the maximality of the matching this neighbour is matched and was added in the fourth round. Together with the edge to its parent, this gives an ear of length two, which we add in the fifth round.

All ears have length at most five, their interiors are new at the time of addition, and the interiors of ears of the same round are disjoint. Consequently, no edge belongs to two ears. When $o$ is replaced by the recorded ends in $B$, the ears remain open or closed ears, since their interiors avoid $B$ and only their two ends may coincide.

Either direction of an ear attaches all its internal vertices strongly to the strong subgraph built before, so the union is strong for every choice of directions. To find a directed path from a vertex towards $B$, we follow its ear in the chosen direction to an end, which is a vertex added in an earlier round or a vertex of $B$, and repeat; for a directed path from $B$ we follow the same rule backwards. At each repetition the round strictly decreases, so there are at most five repetitions, each using at most five edges. This completes the proof.
\end{proof}

We now use the independent directions of the ears to shorten a route whenever its vertices have too large a total degree. The route is always determined by directions which have already been inspected; it need not be a shortest route.

\begin{proposition}\label{prop:density}
Every bridgeless graph $G$ of order $N$ and minimum degree at least $d-1$, where $d\ge 3$, has a strong orientation of diameter at most
\[
\frac{3N}{d}+16000\Bigl(1+\frac{N}{d^2}\Bigr).
\]
\end{proposition}

\begin{proof}
Take $F$ from Lemma~\ref{lem:sparse}, put $e_F=|E(F)|$, and fix a strong orientation of $F$ by Robbins' theorem~\cite{Rob39}. Its diameter is at most $e_F$, also when $F$ is a single vertex. Choose the ears of Lemma~\ref{lem:ears}. We give each ear an independent fair bit which chooses its direction, and we give each remaining edge outside $F$ its own independent fair bit. In particular, edges with both ends in $B=V(F)$ are not fixed unless they belong to $F$. Every assignment of all these bits gives a strong orientation.

Fix an ordered pair of vertices. By inspecting at most ten ear bits, we find its two directed paths to and from $B$ and join their ends by a fixed directed path in $F$. Erasing closed portions gives a directed path of length at most
\[
\ell_0=\frac{18N}{d}+50 ,
\]
all of whose edges have known directions. We fix deterministic rules for all subsequent choices and for the erasing of closed portions. Put
\[
p=\Bigl\lceil 40\bigl(\ell_0+\ln(N+2)+2\bigr)\Bigr\rceil\qquad\text{and}\qquad\Theta=2e_F+10(10+2p).
\]
We shall inspect at most $p$ further candidates, each requiring at most two new bits.

\emph{Finding a fresh candidate.} At any step, we temporarily delete all edges of $F$ and all edges controlled by bits which have already been inspected. An ear bit controls at most five edges, so throughout these $p$ inspections at most $e_F+5(10+2p)$ edges are deleted. For the current directed path $P=(v_0,\dots,v_\ell)$ put
\[
\operatorname{vol}(P)=\sum_{v\in V(P)}\bigl(\deg_G(v)+1\bigr).
\]
Deleting an edge decreases the total size of the closed neighbourhoods of the vertices of $P$ by at most two. We split the vertices of $P$ into three classes according to their indices modulo three. If $\operatorname{vol}(P)>3N+\Theta$, the total size of the closed neighbourhoods in the remaining graph is still greater than $3N$, so some class has total size greater than $N$. Hence two different vertices $v_i,v_j$ of that class, with $i<j$, have intersecting closed neighbourhoods in the remaining graph; in particular $j-i\ge 3$. If the intersection gives a direct edge $v_iv_j$, we inspect its bit and try to use it forwards. Otherwise we choose a common neighbour $z$ and try the route $v_i,z,v_j$. The selected edges were not deleted, so their bits are new. If they have different bits, the conditional probability of success is $1/4$. If they have the same bit, both lie on one ear; at their common vertex $z$ they are consecutive edges of this ear, which includes the first and the last edge of a closed ear, and the two directions of the ear give the two opposite routes through $z$. So the required route succeeds with probability $1/2$. A direct edge also succeeds with probability $1/2$.

The candidate is chosen using only the graph, the current path and the inspected bits. Conditional on this information, all bits not yet inspected are still independent and fair. On success, we replace the subpath of $P$ between $v_i$ and $v_j$ by the route of length one or two and erase closed portions; since $j-i\ge 3$, the length strictly decreases. On failure, we keep the current path and record the inspected values. In both cases, the next path and the directions of all its edges are determined by the inspected bits, so the same conditional argument applies at the next step. Its validity does not require the volume of the path to decrease.

\emph{Stopping all requests in one orientation.} We stop as soon as $\operatorname{vol}(P)\le 3N+\Theta$; until then there is always a fresh candidate. For the purpose of the probability estimate only, we fill any unused inspections after stopping with independent fictitious trials of success probability $1/4$. If $S_p$ is the resulting number of successes, conditional expectation gives
\[
\mathbb{E}\,2^{-S_p}\le(7/8)^p .
\]
If the process has not stopped after $p$ inspections, it has had fewer than $\lceil\ell_0\rceil$ successes, since each success shortens a path whose initial length is at most $\ell_0$, and a path with no edge cannot exceed the volume threshold. Therefore the probability of not stopping is at most
\[
2^{\lceil\ell_0\rceil}(7/8)^p<(N+1)^{-3}
\]
by Markov's inequality; the last inequality follows from $\ln 2<1$, $\ln(8/7)>1/8$ and the definition of $p$. Each request is analysed on the same probability space of edge and ear bits. The inspections for different requests need not be independent, and they do not change the directions. The union bound over at most $N^2$ ordered pairs is less than one, so some assignment of the bits lets every request stop. For each resulting path, every vertex has closed degree at least $d$, so $d(\ell+1)\le 3N+\Theta$, and the diameter is at most $3N/d+\Theta/d$.

Finally,
\[
p\le\frac{720N}{d}+2081+40\ln(N+2)\qquad\text{and}\qquad\Theta\le\frac{14436N}{d}+800\ln(N+2)+41720 .
\]
Using $\ln(N+2)\le 2\sqrt N$, $2\sqrt N/d\le 1+N/d^2$ and $d\ge 3$, we obtain
\[
\frac{\Theta}{d}\le 15236\frac{N}{d^2}+800+\frac{41720}{3}<16000\Bigl(1+\frac{N}{d^2}\Bigr).
\]
This completes the proof.
\end{proof}

For the main argument, the preceding bound must preserve a prescribed core and charge only the vertices outside its closed neighbourhood. The next construction reduces this relative statement to Proposition~\ref{prop:density} without replacing actual vertices of the core by free transfers.

\begin{theorem}\label{thm:density}
Let $G$ be a bridgeless graph of minimum degree at least $d-1$, where $d\ge 3$, and let $K$ be a core of $G$ whose prescribed orientation has diameter $\Delta_K$. Put $q=|V(G)\setminus N_G[K]|$. Then the prescribed orientation extends to a strong orientation of $G$ of diameter at most
\[
\Delta_K+\frac{3q}{d}+33000\Bigl(1+\frac{q}{d^2}\Bigr).
\]
\end{theorem}

\begin{proof}
Write $W_1=N_G(K)\setminus K$ and $W_2=V(G)\setminus N_G[K]$, so that $q=|W_2|$. We first enlarge the prescribed core by a bounded amount in directed distance. Let $\Gamma_W$ be the auxiliary graph on $W_1$ in which two distinct vertices are adjacent if they are adjacent in $G$ or have a common neighbour in $W_2$.

In each component of $\Gamma_W$ with at least two vertices, choose a spanning tree and its two colour classes $X^+$ and $X^-$. For every vertex of the component, choose one edge to $K$, directed away from $K$ at vertices of $X^+$ and towards $K$ at vertices of $X^-$. Realise every edge of the spanning tree either by an edge of $G$ directed from $X^+$ to $X^-$, or by a path $a\to z\to b$ with $a\in X^+$, $b\in X^-$ and $z\in W_2$. These requirements are consistent even if one vertex $z$ is used several times, since all used edges between $X^+$ and $z$ enter $z$ and all used edges between $z$ and $X^-$ leave $z$. Moreover, the neighbours in $W_1$ of a fixed vertex $z\in W_2$ lie in one component of $\Gamma_W$, so different components cannot impose conflicting requirements.

Every vertex of $X^+$ has a directed path from $K$ of length one and a directed path to $K$ of length at most three, and the reverse holds for $X^-$. Every used vertex $z$ has directed paths to and from $K$ of length at most two. An isolated vertex of $\Gamma_W$ with at least two edges to $K$ is also added, with two oppositely directed edges to $K$. Let $K^+$ consist of $K$, all these added vertices of $W_1$, and all used vertices $z$. These directions make $G[K^+]$ strong with diameter at most $\Delta_K+6$, using two boundary routes and at most one route inside $K$. We orient the other edges of $G[K^+]$ arbitrarily, keeping the prescribed directions inside $K$.

Every vertex of $W_1$ outside $K^+$ is an isolated vertex of $\Gamma_W$ with exactly one edge to $K$. It has no neighbour in $W_1$, and hence at least $d-2$ neighbours in $W_2$. These sets of neighbours in $W_2$ are disjoint for different such vertices, since these vertices are isolated in $\Gamma_W$. Hence, with $m=|V(G)\setminus K^+|$,
\[
m\le q+\frac{q}{d-2}.
\]
If $m=0$, the bound already follows.

Contract $K^+$ to a root, keeping all edges. The resulting multigraph is connected and bridgeless, and every other vertex keeps its degree when parallel edges to the root are counted with multiplicity. In each class of $\ell$ parallel edges between a vertex and the root we retain exactly $\min\{\ell,d-1\}$ edges. A vertex which loses edges still has $d-1$ edges to the root, and a changed class keeps at least two edges, so connectivity and bridgelessness are preserved: a cut through a changed class still contains at least two edges, and every other cut is unchanged. Now replace the root by a new complete graph on $d$ vertices, and attach the kept root edges of every vertex to distinct vertices of this complete graph, keeping their labels. This gives a simple graph $G^\#$ of order $N=m+d$ and minimum degree at least $d-1$. It is connected and bridgeless: the edges of the complete graph lie in triangles, and for every other edge $e$, a route avoiding $e$ in the truncated multigraph lifts to $G^\#$ by joining, whenever it passes through the root, the two attachment vertices inside the complete graph.

Apply Proposition~\ref{prop:density} to $G^\#$ and contract the new complete graph back to the root. Projecting directed paths gives a strong orientation of the truncated quotient with no larger diameter; orient the discarded parallel edges arbitrarily. Lemma~\ref{lem:lift}, applied to the core $K^+$ of diameter at most $\Delta_K+6$, now gives a strong orientation of $G$ of diameter at most
\[
\Delta_K+6+\frac{3N}{d}+16000\Bigl(1+\frac{N}{d^2}\Bigr).
\]
Since $N\le q+q/(d-2)+d$, this is at most
\[
\Delta_K+\frac{3q}{d}+\Bigl[9+16000\Bigl(1+\frac1d\Bigr)\Bigr]+\Bigl[\frac{3d}{d-2}+16000\Bigl(1+\frac{1}{d-2}\Bigr)\Bigr]\frac{q}{d^2}.
\]
For $d\ge 3$, both bracketed constants are less than $33000$. This completes the proof.
\end{proof}

\section{Orienting closed trails through one root}\label{sec:trails}

Two edge-disjoint paths from a root to a vertex form a closed trail when one path is followed forwards and the other backwards. Different pairs may share edges, so their chosen directions cannot simply be combined. The following lemma resolves this conflict. Its bound depends logarithmically on the number of trails, and not on the number of vertices which they cover.

\begin{lemma}\label{lem:trails}
Let $\mathcal{U}$ be the union of $k\ge 1$ closed trails, each of which contains a specified vertex $r$ and has length at most a positive integer $L$. Parallel edges are allowed and regarded as distinct. Then $\mathcal{U}$ has a strong orientation of diameter at most
\[
L\bigl(10+6\ln(1+k)\bigr).
\]
\end{lemma}

\begin{proof}
Choose a direction of traversal for each trail, and assign to each trail an independent number uniformly distributed in $(0,1)$. With probability one these numbers are distinct. We orient each physical edge according to the traversal of the trail with the smallest number among the trails containing that edge.

Every assignment of distinct numbers gives a strong orientation. To find a directed walk from a vertex $x$ to $r$, we reveal the numbers of all trails through $x$ and start along the trail with the smallest revealed number. At each subsequent vertex we reveal the numbers of all trails through it which are still unknown. We change to another trail only if a newly revealed number is smaller than the current number, and then we choose the trail with the smallest such number. On changing trails, we use a predetermined occurrence of the current vertex in the new trail; without a change, we continue from the current occurrence in the current trail and do not reset the position when a vertex is repeated. The current trail always has the smallest number revealed so far. Every trail containing the next edge passes through the current vertex, so the numbers of all trails competing for this edge have already been revealed, and the edge has the required direction. If no further change occurs, the current trail reaches $r$ within $L$ steps. Each change strictly decreases the current number, so the process ends at $r$. Following the traversal directions backwards by the same rule gives, read in reverse, a directed walk from $r$ to $x$. This proves strong connectivity.

The number of changes has a uniform probability bound. Consider one of these two searches at a fixed vertex. We order every newly revealed group of trails by fixed labels, and after the search ends we reveal all remaining numbers in the order of their labels. The identity of the next unknown number depends only on the numbers already revealed, and conditional on this information, every unknown number is still uniform and independent of the others. Hence the resulting sequence of $k$ values is a sequence of independent uniform random variables. Let $\chi$ be the number of record minima in this sequence, including the first value. Each followed segment of the search starts at a record minimum, so the search has at most $\chi$ segments, each of length at most $L$. The relative rank of the $j$-th value among the first $j$ values is uniform on $\{1,\dots,j\}$, independently for different $j$, as one sees by inserting the values one by one in their relative order. Consequently,
\[
\mathbb{E}\,2^{\chi}=\prod_{j=1}^{k}\Bigl(1+\frac1j\Bigr)=k+1 .
\]
Hence, for every integer $a\ge 0$, each of the two directed distances between $x$ and $r$ exceeds $aL$ with probability at most $(k+1)2^{-a}$.

A vertex far from $r$ supplies many witnesses. The union $\mathcal{U}$ has at most $kL+1$ vertices. If the distance from some vertex to $r$ exceeds $(a+1)L$, then a shortest directed path from this vertex to $r$ contains at least $L$ distinct vertices whose distances to $r$ exceed $aL$, and the same holds for distances from $r$. Let $Z$ be the number of vertices whose distance to or from $r$ exceeds $aL$, counting the two directions separately. The bounds above give
\[
\mathbb{E}Z\le 2(kL+1)(k+1)2^{-a}.
\]
Thus the probability that some distance to or from $r$ exceeds $(a+1)L$ is at most $\mathbb{E}Z/L\le 2(k+1)^22^{-a}$. Take $a=\lceil\log_2(8(k+1)^2)\rceil$. Then this probability is at most $1/4$, so some assignment makes every distance to and from $r$ at most $(a+1)L$, and its diameter is at most $2(a+1)L$. Finally, $2(a+1)\le 10+4\log_2(k+1)\le 10+6\ln(k+1)$. The searches for different vertices are only analyses of the same assignment; their independence is neither claimed nor needed. This completes the proof.
\end{proof}

We next bound how many trails are needed at a given length. The selection is undirected. Recall that, for a connected vertex set $K$ of a bridgeless graph $H$ and a vertex $x\notin K$, $\mu_K(x)$ denotes the minimum total length of two edge-disjoint paths from the contracted root to $x$ in $H/K$. This number is finite by the edge version of Menger's theorem; the paths need not lift to edge-disjoint paths inside $K$.

\begin{lemma}\label{lem:batch}
Let $H$ be a bridgeless graph of minimum degree at least $d-1$, where $d\ge8$. Let $K$ be a nonempty connected set, and write $Q=H/K$ with root $r$. Suppose that $\mu_K(x)\le L$ for every $x\notin K$, and fix an integer $s\ge100$. Start with $U=\{r\}$. Whenever a vertex $x\notin U$ has minimum two-path total length greater than $s$ in $Q/U$, choose two edge-disjoint $r$--$x$ paths of minimum total length in the fixed graph $Q$, and add their vertices to $U$. Each addition enlarges $N_Q[U]$ by at least $ds/200$ ordinary vertices. Consequently, if $q=|V(H)\setminus N_H[K]|$, then at most $200q/(ds)$ additions are made, and at termination every vertex has two-path total length at most $s$ from the root in $Q/U$.
\end{lemma}

\begin{proof}
Consider one addition, with paths $P_1,P_2$ and support $W$. The support-degree estimate in the proof of Lemma~\ref{lem:balls} gives at most five distinct neighbours in $W$ for any vertex of $W$. The root has at most two: its used neighbours are the first vertices of the two paths, and every unused root--support edge is parallel to a used edge. A vertex outside $W$ has at most six neighbours in $W$, since all its incident edges to $W$ are unused, their other ends have three possible labels, and each label occurs at most twice.

Put $h=\lfloor s/2\rfloor-2$ and read both paths from $x$ towards $r$. If a path meets $U$ within $h$ edges, its initial part is a short route to $U$. Otherwise its first $h+1$ vertices lie outside $U$. A vertex $z\in N_Q[U]\setminus W$ is an \emph{available neighbour} of this initial part if it is adjacent to one of its vertices $u$. If $z\in U$, append $uz$ to the initial part up to $u$; otherwise choose $zt$ with $t\in U$ and append $uz,zt$. Stopping at the first visit of $U$ gives a path of length at most $h+2$ with no internal vertex in $U$.

Suppose that each path either meets $U$ within $h$ edges or has at least two distinct available neighbours. If both need available neighbours, choose different ones. The original initial parts have disjoint edges. Every appended edge has an end outside $W$, so belongs to neither initial part. The first appended edges are distinct because their ends outside $W$ are distinct, and the final appended edges, when present, are distinct because their ends outside $U$ are distinct. A first appended edge of one route cannot be the final appended edge of the other: its end on the initial part lies in $W\setminus U$, whereas the latter edge has an end in $U$ and its other end outside $W$. After contracting $U$, the two routes therefore give edge-disjoint paths of total length at most $2(h+2)\le s$, a contradiction.

Thus one path has $h+1$ initial vertices outside $U$ and at most one available neighbour. Discarding the at most two root neighbours leaves at least $h-1$ vertices. Each has no neighbour in $K$, so it retains at least $d-1$ distinct ordinary neighbours in $Q$. At most five are in $W$ and at most one is in $N_Q[U]\setminus W$, leaving at least $d-7$ neighbours outside $W\cup N_Q[U]$. Each such new vertex is counted at most six times. Since $s\ge100$ and $d\ge8$,
\[
\bigl|N_Q[U\cup W]\setminus N_Q[U]\bigr|
\ge\frac{(h-1)(d-7)}6
\ge\frac{(0.46s)(d/8)}6
>\frac{ds}{200}.
\]
Here $h-1=\lfloor s/2\rfloor-3\ge0.46s$.

The successive new neighbourhood sets are disjoint and lie among the $q$ ordinary vertices initially outside $N_Q[r]$. Hence there are at most $200q/(ds)$ additions. Each adds its chosen $x$, so the process terminates with the stated stopping condition. All pairs are chosen in the fixed graph $Q$, and thus have the same root and total length at most $L$.
\end{proof}

\section{Completion at decreasing lengths}\label{sec:completion}

We use the following radius theorem of Chvátal and Thomassen~\cite{CT78}: if every vertex of a bridgeless graph is within distance $\rho$ of a vertex $v$, then the graph has an orientation in which every vertex has directed distance at most $\rho^2+\rho$ to $v$ and from $v$. Since contracting a core creates parallel edges, we record precisely how we use this theorem.

\begin{corollary}\label{cor:radius}
Let $K$ be a core of a bridgeless graph $H$ whose prescribed orientation has diameter $\Delta_K$. If every vertex of $H$ is within distance $R$ of $K$, where $R\ge0$ is an integer, then the prescribed orientation extends to a strong orientation of $H$ of diameter at most
\[
\Delta_K+2(R+1)(R+2).
\]
\end{corollary}

\begin{proof}
Contract $K$ to a vertex $o$, keeping parallel edges. The resulting multigraph $\Gamma=H/K$ is connected, bridgeless and loopless, and every vertex is within distance $R$ of $o$. From every class of $\ell$ parallel edges we keep exactly $\min\{\ell,2\}$ edges; as in the proof of Theorem~\ref{thm:density}, this preserves connectivity and bridgelessness, since a changed class keeps two edges, and it changes no distance. We replace $o$ by a triangle $o_1o_2o_3$ and, for every neighbour of $o$, attach its kept edges to $o$ to distinct vertices of the triangle; all other edges remain unchanged. The resulting graph $H^\#$ is simple, connected and bridgeless: the edges of the triangle lie in a triangle, and a route avoiding a given edge in the truncated multigraph lifts to $H^\#$ by joining its attachment vertices inside the triangle. Every vertex of $H^\#$ is within distance $R+1$ of $o_1$. By the theorem of Chvátal and Thomassen, $H^\#$ has an orientation in which all directed distances to $o_1$ and from $o_1$ are at most $(R+1)^2+(R+1)=(R+1)(R+2)$.

We contract the triangle back to $o$ and orient the deleted parallel edges arbitrarily. Projecting directed paths, and stopping at the first visit of $o$, respectively starting at the last visit of $o$, shows that in the resulting orientation of $\Gamma$ all directed distances to $o$ and from $o$ are at most $(R+1)(R+2)$. Concatenating a route to $o$ with a route from $o$ shows that this orientation has diameter at most $2(R+1)(R+2)$, and Lemma~\ref{lem:lift} gives the asserted bound.
\end{proof}

The next estimate completes the orientation when many vertices remain outside the closed neighbourhood of the core.

\begin{proposition}\label{prop:decreasing}
Let $H$ be a bridgeless graph of minimum degree at least $d-1$, where $d\ge8$. Let $K_0$ be a core whose prescribed orientation has diameter $\Delta_0$, and suppose that $\mu_{K_0}(x)\le36\,000\,000d$ for every $x\notin K_0$. Put $q_0=|V(H)\setminus N_H[K_0]|$. Then the prescribed orientation extends to a strong orientation of $H$ of diameter at most
\[
\Delta_0+2\cdot10^9d\left[1+\ln\left(1+\frac{q_0}{d^2}\right)\right].
\]
\end{proposition}

\begin{proof}
Put $L_0=36\,000\,000d$. At stage $i$, the current core $K_i$ satisfies $\mu_{K_i}(x)\le L_i$ for every vertex outside it. As long as $s_i=\lfloor L_i/2\rfloor\ge100$, apply Lemma~\ref{lem:batch} in the fixed graph $H/K_i$ with $s=s_i$, choosing the entire batch before assigning any direction. Suppose the batch contains $k_i$ pairs. Each pair forms a closed trail through the root of length at most $L_i$.

If $k_i>0$, apply Lemma~\ref{lem:trails} to their union, and orient the other edges of the quotient induced by the enlarged vertex set arbitrarily; this keeps the orientation strong and does not increase any distance. By Lemma~\ref{lem:lift}, this orientation lifts to the enlarged core, keeping the prescribed core directions, and the diameter increases by at most $L_i[10+6\ln(1+k_i)]$. Contracting the enlarged core gives exactly the graph $Q/U$ in Lemma~\ref{lem:batch}, so the new two-path totals are at most $s_i$; put $L_{i+1}=s_i$. An empty batch requires no enlargement or diameter cost, and we simply use the smaller bound. All edges outside the current core remain unoriented until needed.

By Lemma~\ref{lem:batch}, $k_i\le200q_0/(ds_i)$. Put $t=q_0/d^2$. Repeated halving gives $L_i=\lfloor L_0/2^i\rfloor$ and, at every active stage, $s_i\ge L_0/(4\cdot2^i)$. Thus
\[
k_i\le\frac{800q_0\,2^i}{dL_0}\le t2^i,
\qquad
\ln(1+k_i)\le\ln(1+t)+i\ln2.
\]
Using $L_i\le L_0/2^i$ and $\sum_{i\ge0}2^{-i}=\sum_{i\ge0}i2^{-i}=2$, the total increase over all active stages is at most
\[
L_0[20+12\ln(1+t)+12\ln2].
\]
Here $q_0$ bounds the number of trails at each stage; it is not charged as a fresh volume at each stage.

At termination $L_i<200$, so the shorter path of each pair has length at most $99$. Corollary~\ref{cor:radius} adds at most $2\cdot100\cdot101=20200$. Since $d\ge8$ and $L_0=36\,000\,000d$,
\[
L_0[20+12\ln(1+t)+12\ln2]+20200
\le2\cdot10^9d[1+\ln(1+t)].
\]
Indeed, the coefficient of $\ln(1+t)$ on the left after division by $d$ is $432\,000\,000$, and its constant term is less than $1.02\cdot10^9$. This proves the proposition.
\end{proof}

\begin{proof}[Proof of Theorem~\ref{thm:main}]
If $n=1$, the result is immediate. Otherwise bridgelessness implies $\delta\ge2$. For $\delta=2$, any strong orientation has diameter at most $n-1<3n/(\delta+1)$. Assume henceforth that $\delta\ge3$ and put $d=\delta+1$.

Apply Lemma~\ref{lem:kernel} and let $H$ be the reduced graph, of order $n_H$. Suppose first that $\delta\ge7$, so $d\ge8$. Proposition~\ref{prop:initial} gives a core $K_0$ of prescribed diameter $\Delta_0\le3|N_H[K_0]|/d$ and two-path bound $36\,000\,000d$. Write
\[
q_0=n_H-|N_H[K_0]|,
\qquad t=q_0/d^2.
\]
If $t\ge10^{12}$, use Proposition~\ref{prop:decreasing}. Its extra cost is at most $3q_0/d$: at $t=10^{12}$,
\[
2\cdot10^9[1+\ln(1+t)]<2\cdot10^9\cdot29<3t,
\]
and $t/[1+\ln(1+t)]$ is increasing for $t>0$. The resulting diameter is at most $3n_H/d$. If $t<10^{12}$, instead apply Theorem~\ref{thm:density} to the same core, obtaining diameter at most
\[
\Delta_0+\frac{3q_0}{d}+33000(1+t)
\le\frac{3n_H}{d}+33000(1+10^{12}).
\]
These are alternative constructions selected by the single graph parameter $q_0$, so each supplies one orientation for every ordered pair.

For $4\le\delta\le6$, Proposition~\ref{prop:initial} gives a core $K$ of diameter at most $3|N_H[K]|/d$ from which every vertex is within distance
\[
\frac{36\,000\,000d}{2}\le126\,000\,000.
\]
For $\delta=3$, Corollary~\ref{prop:initialsmall} gives such a core with radius at most $82\,320\,000$. In either case Corollary~\ref{cor:radius} yields diameter at most
\[
\frac{3n_H}{d}+2(R_*+1)(R_*+2),
\qquad R_*=126\,000\,000.
\]

Finally recover an orientation of $G$ by Lemma~\ref{lem:kernel}. Its order term combines with $3n_H/d$ to give $3n/d$, and its additive cost is $26$. Both possible errors satisfy
\[
33000(1+10^{12})+26
=33\,000\,000\,000\,033\,026<3.51\cdot10^{16},
\]
\[
2(R_*+1)(R_*+2)+26
=31\,752\,000\,756\,000\,030<3.51\cdot10^{16}.
\]
This proves the stronger bound claimed after Theorem~\ref{thm:main}, and hence the theorem.
\end{proof}

\section*{Artificial intelligence use declaration}

Anthropic's Claude and OpenAI's Codex were used to assist with computational checks of the explicit constants, the restructuring of several proofs, consistency checks, and the drafting and language editing of this manuscript. The authors take full responsibility for the results presented in this manuscript.


\begin{thebibliography}{9}

\bibitem{AMO93}
R.~K. Ahuja, T.~L. Magnanti and J.~B. Orlin, \emph{Network flows: theory, algorithms, and applications}, Prentice Hall, Englewood Cliffs, NJ, 1993.

\bibitem{BD15}
S.~Bau and P.~Dankelmann, Diameter of orientations of graphs with given minimum degree, \emph{European Journal of Combinatorics} 49 (2015): 126--133.

\bibitem{CT78}
V.~Chvátal and C.~Thomassen, Distances in orientations of graphs, \emph{Journal of Combinatorial Theory, Series B} 24(1) (1978): 61--75.

\bibitem{Coc24}
G.~Cochran, Large girth and small oriented diameter graphs, \emph{Discrete Mathematics} 347(4) (2024): 113846.

\bibitem{CW26}
G.~Cochran and Z.~Wang, On the oriented diameter of graphs with given minimum degree, \emph{Journal of Graph Theory}, Early View (2026). \url{https://doi.org/10.1002/jgt.70099}.

\bibitem{EPPT89}
P.~Erd\H{o}s, J.~Pach, R.~Pollack and Zs.~Tuza, Radius, diameter, and minimum degree, \emph{Journal of Combinatorial Theory, Series B} 47(1) (1989): 73--79.

\bibitem{Rob39}
H.~E. Robbins, A theorem on graphs, with an application to a problem of traffic control, \emph{The American Mathematical Monthly} 46(5) (1939): 281--283.

\bibitem{Sur17}
M.~Surmacs, Improved bound on the oriented diameter of graphs with given minimum degree, \emph{European Journal of Combinatorics} 59 (2017): 187--191.

\bibitem{ST84}
J.~W. Suurballe and R.~E. Tarjan, A quick method for finding shortest pairs of disjoint paths, \emph{Networks} 14(2) (1984): 325--336.

\end{thebibliography}
\end{document}